\documentclass[psamsfonts,a4paper]{amsart}

\usepackage{amssymb,amscd,amsmath,amsthm,color}

\chardef\bslash=`\\ 

\newtheorem{theorem}{Theorem}[section]
\newtheorem{corollary}[theorem]{Corollary}
\newtheorem{lemma}[theorem]{Lemma}
\newtheorem{proposition}[theorem]{Proposition}

\theoremstyle{remark}
\newtheorem{remark}[theorem]{Remark}

\theoremstyle{definition}
\newtheorem{definition}[theorem]{Definition}

\numberwithin{equation}{section}

\newcommand{\thmref}[1]{Theorem~\ref{#1}}

\newcommand{\proref}[1]{Proposition~\ref{#1}}
\newcommand{\lemref}[1]{Lemma~\ref{#1}}
\newcommand{\corref}[1]{Corollary~\ref{#1}}
\newcommand{\remref}[1]{Remark~\ref{#1}}

\newcommand{\N}{\mathbb N}
\newcommand{\Z}{\mathbb Z}
\newcommand{\Zp}{\mathbb Z_p}
\newcommand{\Q}{\mathbb Q}
\newcommand{\R}{\mathbb R}
\newcommand{\C}{\mathbb C}

\newcommand{\inv}{^{-1}}
\newcommand{\id}{\operatorname {id}}
\newcommand{\kk}{\mathcal K}
\newcommand{\rr}{\mathcal R}
\newcommand{\ks}{\kk^*}
\newcommand{\oo}{\mathcal O}
\newcommand{\ox}{\oo^\times}	
\newcommand{\os}{\oo^*}
\newcommand{\osbar}{\overline{\os}}
\newcommand{\ov}{\oo_v}
\newcommand{\kmo}{\kk/\oo}
\newcommand{\mk}{M_\kk}
\newcommand{\mko}{M_\kk^0}

\newcommand{\heck}{C^*_r(\pk,\po)}
\newcommand{\hcal}{\mathcal H}
\newcommand{\teta}{\theta}
\newcommand{\alphin}{\beta}

\newcommand{\po}{P_\oo}
\newcommand{\pk}{P_\kk}

\newcommand{\lspa}{\operatorname{span}}

\def\E{\mathrm e}
\def\ER{\mathrm f}

\newcommand{\mathfrac}{\frac}

\newcommand{\matr}[2]{\left(\begin{matrix}1&#1 \\  
0&#2\end{matrix}\right)}
\newcommand{\smatr}[2]{\bigl(\mbox{\tiny{$\begin{matrix}1&#1 \\  
0&#2\end{matrix}$}}\bigr)}
\newcommand{\stabo}[1]{\os_{\!#1}}

\newcommand{\divides}{\operatorname*{|}}
\newcommand{\aid}{\mathfrak{a}}	
\newcommand{\bid}{\mathfrak{b}}	

\def\diff{\mathfrak{D}}
\def\invd{\mathcal{D}\inv}
\def\invdv{\mathfrak{D}\inv_v}

\newcommand\sgn{\operatorname*{sgn}}

\newcommand{\pid}{\mathfrak{p}}
\newcommand\Cl{\mathrm{Cl}}
\newcommand{\convolve}{\operatorname*{\ast}}
\newcommand{\Tr}{\operatorname{Tr}}
\newcommand{\supp}{\operatorname{supp}}

\newcommand{\Gal}{\operatorname{Gal}}
\newcommand\A{\mathbb{A}}	

\begin{document}

\title[Phase transitions on Hecke C*-algebras]{Phase transitions on  
Hecke C*-algebras\\ and class-field theory over $\Q$}

\date{September 28, 2004.}
\author{Marcelo Laca}
\thanks{Supported by the National Science and Engineering  
Research Council of Canada}

\address{Department of Mathematics and Statistics,
        University of Victoria, Victoria, BC, V8W 3P4, CANADA}
\email{laca@math.uvic.ca}

\author{Machiel van Frankenhuijsen}
\address{Department of Mathematics,
Utah Valley State College,
Orem, UT 84058-5999, USA}
\email{vanframa@uvsc.edu}


\begin{abstract}
We associate a canonical Hecke pair of semidirect product groups
to the ring inclusion of the algebraic integers $\oo$
in a number field $\kk$, and we construct a C*-dynamical system on the  
corresponding Hecke
C*-algebra, analogous to the one constructed by Bost and Connes for the
inclusion of the integers in the rational numbers.
We describe the structure of the resulting Hecke C*-algebra as a  
semigroup crossed product
and then, in the case of class number one,
analyze the equilibrium (KMS) states of the dynamical system.
The extreme KMS$_\beta$ states
at low-temperature exhibit a phase transition with symmetry breaking
that strongly suggests a connection with class field theory. Indeed,
for purely imaginary fields of class number one,
the group of symmetries, which acts freely and transitively
on the extreme KMS$_\infty$ states, is isomorphic to
the Galois group of the maximal abelian extension over the field.
However, the Galois action
on the restrictions of extreme KMS$_\infty$ states to the (arithmetic)  
Hecke algebra over $\kk$,
as given by class-field theory, corresponds to the action of the  
symmetry group
if and only if the number field $\kk$ is $\Q$.
\end{abstract}

\maketitle

\section*{Introduction}
The C*-dynamical system based on a noncommutative Hecke C*-algebra  
constructed by Bost and  Connes
from the inclusion of the integers in the rationals has inspired  
several authors to construct Hecke C*-algebras
associated to algebraic number fields and function fields, e.g.  
\cite{har-lei,alr,coh}.
These constructions share many of the interesting features of the  
Bost-Connes
construction, in particular they all have a phase transition with  
spontaneous symmetry breaking
at low temperature, a partition function related to the zeta function  
of the number field,
and a unique type III factor equilibrium state at high temperature.
   The Bost-Connes C*-dynamical system has a group of symmetries that  
acts freely and transitively on the extreme KMS$_\beta$ states for  
$\beta>1$, and which
   is isomorphic to the Galois group of  the maximal abelian extension $\Q^{ab}$  of  
$\Q$.
Moreover,  the embeddings of $\Q^{ab}$ in $\C$ are
given by the evaluation of the extreme KMS$_\infty$ states on an  
`arithmetic'
Hecke algebra over $\Q$, and thus their model also has an arithmetic  
symmetry, in which the
Galois group acts by Galois automorphisms on the values of extreme  
KMS$_\infty$ states.
This feature of having concrete Galois groups as symmetries has not  
been generalized to number fields, although in some cases the broken  
symmetries do have an interpretation as abstract Galois symmetries,  
e.g. \cite{har-lei,coh}.

In this work we were initially motivated by the observation that the  
symmetry
groups of \cite{har-lei} are
isomorphic to Galois groups for the nine imaginary quadratic fields of  
class number $1$,
see Remark \ref{lastremark}.
  We study the Hecke C*-algebra
canonically associated to the inclusion of the ring of integers $\oo$ in
an algebraic number field $\kk$. Specifically, we consider the full `$ax+b$'  
group $\oo \rtimes \os$
  of the ring of algebraic integers in that of the field,
$\kk \rtimes \ks$. By `canonical' here we mean that no cross section
is chosen to deal with the presence of the units $\os $ in $\oo$.
Although the inclusion of the full group of units in the subgroup
causes some technical difficulties with the Hecke inclusion
at the onset, these are offset later by an eventual simplification
in the computation of the equilibrium states. Indeed,
the unit group, which when infinite presented an obstacle
to the computation of KMS states on the Hecke algebras of \cite{alr}, is
  `factored out' from our Hecke algebra, and as a result, it is the  
semigroup
of principal integral ideals that acts as
our renormalization semigroup of endomorphisms.

In Section \ref{S: Hecke} we carry out the analysis of the structure of  
our Hecke
algebra for a general number field $\kk$.
We give a presentation in terms of generators and relations and also a  
realization
as the crossed product of a certain commutative *-algebra
by an action of the semigroup of principal integral ideals of $\kk$.
This commutative  *-algebra can be described in three different ways:  
first
as the fixed point algebra of the action
of the units on the group algebra of $\kk/\oo$, then
as the Hecke *-algebra of an intermediate Hecke pair, and finally, via  
the
Fourier transform in Section 2, as an algebra
of continuous functions on the dual of $\kmo$, which we view
as the `adelic global inverse different'.
Similar descriptions of the corresponding C*-algebra are also possible
thanks to a natural profinite compactification of the group $\os$ of units
that arises from their action on $\kmo$,  which allows us to consider
the appropriate compact orbit space, see
\lemref{aside}.

The Hecke C*-algebra has a natural dynamics and  a natural group of  
`geometric' symmetries.
These geometric symmetries are induced from the symmetries of the space  
of orbits of the additive classes in $\kmo$ under the multiplicative  
action of the unit group $\os$.
To describe the phase transition of KMS states, in Section  
\ref{S:classnumberone},
we restrict our attention to fields of class number one.
The phase transition and the symmetries for fields of higher class  
numbers
require techniques beyond the ones developed in \cite{diri} for lattice  
semigroups,
and will be considered elsewhere.

Our main result is \thmref{main}, where we show that for fields of  
class number one
there is a unique KMS$_\beta$ state for $0\leq \beta\leq 1$,
while for $\beta > 1$ the symmetry group acts freely and transitively
on extreme KMS$_\beta$ states.
In the absence of real embeddings, that is, for purely imaginary
fields, this symmetry group is isomorphic to the
Galois group of the maximal abelian extension of the field, but
in general, the symmetry group is missing the $\{\pm 1\}$
factors in the Galois group corresponding to complex conjugation on  
each real embedding.

To explore the possibility of a Galois action of the symmetry group,  
that is,
the possibility of arithmetic symmetries,
we consider in Section \ref{S:Classfields}
 the natural candidate for an arithmetic Hecke algebra over  
the field.
When the extreme KMS$_\infty$ states are evaluated on this algebra, the  
resulting
  values are algebraic numbers, but the Hecke algebra model is
  based on torsion alone, in the sense that
  these numbers always lie in the maximal
{\em cyclotomic} extension, and hence are not enough to support
a free transitive action of the Galois group of the maximal abelian  
extension. Indeed, in \thmref{badnews} when we compare the action of an idele
as a geometric symmetry to its Galois action (via the Artin map) as an  
arithmetic symmetry on KMS$_\infty$ states,
  we discover that they are intertwined only when $\kk$ is equal $\Q$.

This leads us to the conclusion that the connection with class field  
theory
in the Bost--Connes Hecke C*-algebra model is an exceptional feature
  derived from the peculiarities of the base field $\Q$;
  specifically, a consequence of the Kronecker-Weber Theorem,
by which the maximal abelian extension of $\Q$ is in fact the maximal  
cyclotomic extension.
However, the results of \cite{bos-con} are too suggestive and the  
possibility of
a noncommutative model that supports both arithmetic and geometric
symmetries for number fields other than $\Q$
is too tempting to be abandoned without further effort. Indeed,
as this work was being finished, we learned of a new construction due
to Connes and Marcolli \cite{con-mar} in which a new
C*-dynamical system is
proposed that promises to generate enough algebraic values for  
KMS$_\infty$ states
to support free transitive Galois actions.

Acknowledgments:
This research was initiated during a visit of M.L.\ to IHES, and carried  
out through
several visits, of M.v.F.\ 
to the University of Newcastle, Australia
  and to the SFB 478 at the University of Muenster,
  and of M.L.\ to Utah Valley State College.
  We acknowledge the support of those institutions and we would like to  
thank the
  respective mathematics departments for their hospitality.

\section{The Hecke C*-algebra}
\label{S: Hecke}
We use the following notation: $\kk$ will denote an algebraic number  
field with ring of integers $\oo$.
The invertible elements of $\kk$ form a multiplicative group, which  
will be denoted by $\ks$;
since $\kk$ is a field,  these are simply the nonzero elements of $\kk$.
Similarly, the group of invertible elements (units) of $\oo$ will be  
denoted by $\os$;
notice that $\os$ is strictly smaller than the multiplicative semigroup  
$\ox$ of  nonzero integers.

To the inclusion of $\oo$ in $\kk$  we associate the following  
inclusion of `$ax+b$ groups'
\[
\po := \matr{\oo}{\os} \ \subset \
\pk := \matr{\kk}{\ks}  .
\]
We shall prove that this is a Hecke (or almost normal)  inclusion; in  
other words,
that every double coset
contains finitely many right cosets.
This is equivalent to saying  that the action of $\po$
by right-multiplication on the right coset space $\po\backslash \pk$  
has finite orbits.
Before we can prove this, we need to develop some basic notation and  
results
concerning the multiplicative action of the units $\os$ on $\kk$ modulo  
$\oo$.
For each $r\in \kk$,
denote by $\stabo{r}$ the subgroup of $\os$ that fixes $r$ modulo $\oo$:
\begin{gather}
\stabo{r}=\{ u\in\os\colon ru=r\bmod\oo\} .
\end{gather}
More generally, for each fractional ideal $\aid$ of $\kk$, denote by  
$\stabo{\aid}$
the subgroup of $\os$ that fixes each element of $\aid$ modulo $\oo$:
\begin{gather}
\stabo{\aid} :=\{ u\in\os\colon ru = r\bmod\oo\text{ for each  
}r\in\aid\} =\bigcap_{r\in\aid}\stabo{r}.
\end{gather}

In particular, note that $\stabo{r}=\stabo{r\oo}$, that
$\stabo{r}$ depends only on the class of $r$ in $\kk/\oo$,
and  that $\stabo{\aid}$ depends only on the ideal $\aid +\oo$.
Also, if $\aid$ and $\bid$ are ideals of $\kk$ with
$\aid\subset\bid$, then $\stabo{\aid}\supset\stabo{\bid}$.

\begin{lemma}
For every fractional ideal\/ $\aid$ in $\kk$, the index of\/  
$\stabo{\aid}$ in $\os$ is finite.
\end{lemma}

\begin{proof}
Let  $\bid:= (\aid +\oo )\inv$. Then $\bid$ is an integral ideal and  
the multiplicative group
$(\oo / \bid)^*$ is finite. Let $n$ be the exponent of $(\oo /  
\bid)^*$, so that
  $u^n = 1\bmod\bid$ for every $u\in\os$.
It follows that $ru^n = r\bmod\oo$ for every $r\in \aid$,
hence $\os / \stabo{\aid}$ has exponent at most $n$.
Since~$\os$ is finitely generated, this shows that the quotient $\os /  
\stabo{\aid}$
is finite.
\end{proof}

The right coset of $\gamma := \smatr{y}{x}$ is
\begin{gather*}
\po\gamma = \left\{\matr{y+xa}{xu}\colon a\in\oo, u\in\os\right\}
=\matr{y+x\oo}{x\os} ,
\end{gather*}
with the obvious notation for a set of matrices.
The image of $\po\gamma$ under the action of
$\smatr{a}{u} \in \po$ by right-multiplication is simply
$
\po\smatr{y}{x} \smatr{a}{u} = \smatr{a+yu+x\oo}{x \os}
$
  and the corresponding right-orbit is
\begin{gather*}
\left\{\matr{a+yu+x\oo}{x \os}  \colon a\in \oo, u\in \os\right\}.
\end{gather*}

Next we compute explicitly the number $R(\gamma)$ of right cosets in the
double coset $\po \gamma \po$ for each $\gamma \in \pk$.

\begin{lemma}\label{L:almostnormal}
Suppose $\gamma = \smatr{y}{x} \in \pk$, and let $\aid := (\oo +  
x\oo)\inv$ be the inverse
of the fractional ideal $\oo + x\oo$.
Then
\begin{gather*}
R(\gamma) = [\os : \stabo{y \aid}] \cdot [\oo : (\oo \cap x\oo)].
\end{gather*}
In particular, the pair $(\pk , \po)$ is a Hecke pair.
\end{lemma}

\begin{proof}
Suppose $u$ and $v$ are in $\os$.
It is easy to see that $u $ and $v$ are in the same class modulo  
$\stabo{y/x}$,
if and only if, for any $a\in\oo$, we have that $a+yu = a+yv\bmod x\oo$.
Hence, fixing $a$, we see that
\[
\matr{a+yu+x\oo}{x\os} \neq \matr{a+yv+x\oo}{x\os}
\]
  if and only if $u\stabo{y/x} \neq v\stabo{y/x}$.
On the other hand, fixing $u\in\os / \stabo{y/x}$ we see that
\[
\matr{a+yu+x\oo}{x\os} \neq
\matr{b+yu+x\oo}{x\os}
\]
  only if $a - b \notin x\oo$, that is,  only if $a \neq b \bmod (\oo  
\cap x\oo)$.
  It follows that every right coset in the right-orbit of $\po \gamma$
arises from a group element $\smatr{a}{u}$ with
$u$ ranging over a set of representatives of $\os / \stabo{y/x}$ and
$a$ ranging over a set of representatives of $\oo / (\oo\cap x\oo)$.
Hence $R(\gamma) \leq |\os / \stabo{y/x}| \cdot |\oo / (\oo\cap x\oo)|  
< \infty$,
proving that $(\pk, \po)$ is a Hecke pair.

This argument only gives an upper bound for $R(\gamma)$ because
the same right coset may arise from different combinations of
representatives of $\os / \stabo{y/x}$ and of $\oo / (\oo\cap x\oo)$.
To account for the redundancy, and to compute $R(\gamma)$,
let  $u_1,u_2\in\os /\stabo{y/x}$ and $a_1,a_2\in\oo /(\oo\cap x\oo )$.
The right cosets corresponding to $\smatr{a_1}{u_1}$  and
$\smatr{a_2}{u_2}$ coincide if and only if
$a_1 + yu_1 + x \oo = a_2 + yu_2 + x \oo $, equivalently, if and only  
if $ a+y(u-1)\in x \oo$,
where we have put $u := u_2 / u_1$ and  $a := (a_2 - a_1)/u_1$.
Note that for each $u$, this is possible for at most one value of $a$  
in $\oo /(\oo\cap x\oo )$.
Furthermore,
\begin{align*}
  a+y(u-1)\in x \oo \quad \text{ for some } a & \\
   \iff \quad  &y (u-1) \in \oo +x \oo    \\
   \iff \quad & z y(u-1) \in \oo  \qquad \text{ for every }  z\in \aid   
\\
   \iff \quad & (zy) u = zy \bmod \oo \quad \text{ for every }  z\in  
\aid  \\
   \iff \quad & u \in \stabo{y\aid}.
\end{align*}
Thus for each fixed pair of representatives $(u_1,a_1)$ of $(\os /  
\stabo{y/x} ) \times (\oo / (\oo\cap x\oo))$
  there are exactly $|\stabo{y\aid}/\stabo{y/x}|$ pairs $(u_2,a_2)$,
with $u_2 = u u_1$ and $a_2 = a_1 + a u_1$
(where $a$ is the unique element modulo $\oo\cap x\oo$ such that  
$a+y(u-1)\in x\oo$),
such that  $\smatr{a_2}{u_2} $ is in the right coset of  $  
\smatr{a_1}{u_1} $.
Since $[\os :\stabo{y/x}]/[\stabo{y\aid}:\stabo{y/x}]=[\os  
:\stabo{y\aid}]$,
the lemma follows.
\end{proof}

The {\em Hecke algebra} $\hcal(\pk,\po)$ is, by definition, the  
convolution *-algebra of complex--valued
bi-invariant functions on $\pk$ that are supported on finitely many  
double cosets.
The convolution product is defined by
\begin{gather}\label{E: convolutiondef}
f\convolve g(\gamma )=\sum_{\gamma_1\in \po\backslash  
\pk}f(\gamma\gamma_1^{-1})g(\gamma_1),
\end{gather}
which is a finite sum because each double coset contains finitely many  
right cosets;
the involution is given by $f^*(\gamma) = \overline{f(\gamma\inv)}$,  
and the
identity is the characteristic function of $\po$.
The same convolution formula, with $g$ replaced by a square integrable  
function $\xi$
on the space $\po\backslash \pk$
  defines an involutive representation $\lambda$ of $\hcal(\pk,\po)$ on  
the Hilbert space
$\ell^2(\po\backslash \pk)$:
\begin{gather}\label{E: lambdadef}
\lambda(f) \xi = f\convolve \xi \qquad \text{ for }f\in \hcal \text{  
and  } \xi \in \ell^2(\po\backslash \pk).
\end{gather}
We shall refer to $\lambda$ as the {\em Hecke representation} of  
$\hcal(\pk,\po)$;  the norm closure of its image
  $\lambda (\hcal(\pk,\po)) $ is, by definition, the {\em Hecke  
C*-algebra}, and is
  denoted by $\heck$.
  Both $\hcal(\pk,\po)$ and $\heck$ come with a canonical {\em time  
evolution} by
  automorphisms $\{ \sigma_t: t\in \R\}$, given by
  \[
  \sigma_t(f) (\gamma) ={\textstyle  
(\mathfrac{R(\gamma)}{L(\gamma)})^{it}} f(\gamma),
  \qquad \gamma \in \pk, \  t\in \R,
  \]
  where $L(\gamma)$ denotes the number of left cosets in the double  
coset of $\gamma$, and is equal to $R(\gamma\inv)$. A routine  
calculation shows that
  this action is spatially implemented on $\ell^2(\po\backslash \pk)$,
  by the unitary group defined by $(U_t\xi)(\gamma) =   
(\mathfrac{R(\gamma)}{L(\gamma)})^{it} \xi(\gamma)$ for $\xi \in  
\ell^2(\po\backslash \pk)$.
See \cite{kri,bin,bos-con} for details.

To simplify the notation, we often write products in the Hecke algebra  
simply as
$fg$ instead of $f\convolve g$, provided there is no risk of confusion,
 and we
use square brackets to indicate the characteristic function of a subset  
of $\pk$.
Thus, the Hecke *-algebra $\hcal(\pk,\po)$ is the linear span of the  
characteristic functions $[\po \gamma \po]$ of double cosets of  
elements $\gamma $ in  $\pk$.

Next we consider two maps $\mu:\ox \to \hcal(\pk,\po)$ and $\teta:\kk  
\to \hcal(\pk,\po)$ defined as follows.
Let $N_a = |\oo /a\oo|$ be the absolute norm of $a\in \ox$, and let
\begin{gather}\label{E: mua}
\mu_a :=\mathfrac{1}{\sqrt{N_a}}\left[\po\matr{0}{a}\po\right];
\end{gather} for each $r\in \kk$ let
\begin{gather}\label{E: er}
\teta_r :=\mathfrac{1}{R(r)}\left[\po\matr{r}{1} \po\right] ,
\end{gather}
where we use the shorthand notation $R(r)$ for $R(\smatr{r}{1})$, the
number of right cosets in the double coset  $\po \smatr{r}{1}\po $.
  Recall that $R(r) = [\os :\os_r]$  by \lemref{L:almostnormal}.

  \begin{remark}\label{R: factorthru}
Since $\po \smatr{0}{a}\po =  \smatr{\oo}{a\os}$ and $N_{ua} = N_a$ for  
every $u\in \os$,
it is clear that the map $\mu$ factors through the quotient $\ox \to  
\ox/\os$,
and hence can be viewed as an injective map from the semigroup of  
principal
integral ideals into $\hcal(\pk,\po) $. Similarly, the double coset of  
$ \smatr{r}{1}$ is
  \[ \po \matr{r}{1}\po =
\matr{\oo+r\os}{\os} = \matr{(\oo+r)\os}{\os},
  \]
  where no link is implied between the two occurrences of $\os$ in the  
last expression.
  Hence, $\teta$ factors through the
quotient $\kk \to \kmo$. In fact, $\teta_r$ depends on $r$ only through  
the
orbit of $r+\oo\in \kmo$ under the multiplicative action of $\os$ on  
$\kmo$
  We shall denote the set of such orbits by $(\kmo)^{\os}$.
  It is clear that different orbits give different $\teta_r$'s,
  because the corresponding supports are disjoint.
\end{remark}

Next,
we show that our Hecke algebra is universal for its relations,
see \thmref{T: universal presentation}.
We need two lemmas to understand the algebra generated by $\teta_r$,
$r\in\kmo$.
Denote by $\E_r$  the indicator function
\[
\E_r :=\left[\po \matr{r}{1}\right]
\]
  of the right  coset of\/ $\smatr{r}{1}$,
and denote by $\ER_r$ the indicator function
\[
\ER_r: =\left[\matr{r}{1}\po\right]
\]
  of the left coset of\/ $\smatr{r}{1}$.

\begin{lemma}\label{L: er is sum}
  For any subgroup $\mathcal S$ of\/ $\stabo{r}$ of finite index in  
$\os$ we have
\begin{gather}\label{E: er is lefts}
\teta_r =\frac{1}{[\os :\mathcal S]}\sum_{u\in\os \!/\mathcal S}\E_{ur}
=\frac{1}{[\os :\mathcal S]}\sum_{u\in\os \!/\mathcal S}\ER_{ur} .
\end{gather}

\end{lemma}

\begin{proof}
The right orbit of $\po \smatr{r}{1}$ consists of the right cosets
$ \smatr{ru+\oo}{\os} =  \po \smatr{ru}{1} $ for $u\in\os$, and it is  
clear that the
union of these right cosets is the common support of
  both $\teta_r$ and $\sum_{u\in\os \!/\mathcal S}\E_{ur}$.
Further, for $\gamma$ in this support,
there is a unique $u\in\os$ modulo $\stabo{r}$ such that
$\po\gamma =\po \smatr{ru}{1}$,
and, modulo $\mathcal S$,
there are $[\stabo{r}:\mathcal S]$ such $u$'s.
Hence the value of the right-hand side of~\eqref{E: er is lefts} at  
$\gamma$ is
$[\stabo{r}:\mathcal S]/[\os :\mathcal S]=1/[\os :\stabo{r}]$,
which, by definition, is the value of $\teta_r$ at $\gamma$.
The proof of the second equality is entirely analogous.
\end{proof}

The functions $\E_r$  and $\ER_r$ are not bi-invariant, so they are  
not in $\hcal$, but
  $\E_r$ is left--invariant and $\ER_r$ is right--invariant, and it  
still makes sense to compute
the convolution product $\ER_r\convolve\E_s$ using formula~(\ref{E:  
convolutiondef}).
We notice however that this product may fail to be left or  
right--invariant.

\begin{lemma}\label{L: convolution of e and f}
For $r,s\in \kk$ the convolution \/ $\ER_r  \convolve \E_s $  is the  
indicator function
\[\ER_r  \convolve \E_s  =
  \left[\bigcup\nolimits_{u\in\os}\smatr{ru+s+\oo}{u}\right] =
\sum\nolimits_{u\in\os}\left[\smatr{ru+s+\oo}{u}\right].
\]
\end{lemma}

\begin{proof}
Since $\E_s (\gamma_1 )=1$ if and only if the right coset $\gamma_1 $  
is precisely $\po\smatr{s}1$,
we have $\ER_r\convolve\E_s(\gamma )=1$ or $0$ according to whether  
$\ER_r(\gamma\gamma_1^{-1})=1$ or $0$ for
$\gamma_1 =\po\smatr{s}1$.
Now $\ER_r(\gamma\gamma_1^{-1})=1$ if and only if  
$\gamma\gamma_1^{-1}=\smatr{r}1\po$
as a left coset,
which means $\gamma  
=\smatr{r}{1}\smatr{a}{u}\smatr{s}{1}=\smatr{ru+s+a}{u}$ for some
$a\in\oo$ and $u\in\os$. Finally, the characteristic function of the  
union can be replaced by a pointwise sum because the
sets are mutually disjoint, so that only one term is nonzero at each  
point.
\end{proof}

\begin{proposition}\label{P:Ateta characterization}
\begin{enumerate}
\item The *-subalgebra $\mathfrak A$ of\/ $\hcal(\pk,\po)$ generated by  
the set\/
  $\{\teta_r\colon r\in (\kmo)^{\os}\}$ is universal for the relations:
\begin{align*}
&\teta_0 = 1  \\
&\teta_{wr} = \teta_r = \teta_r^*, \qquad r\in \kmo, w\in \os,\\
&\teta_r \teta_s = \mathfrac{1}{R(r)}
\mathfrac{1}{R(s)}\sum_{u\in\os \!/\stabo{r}}\sum_{v\in\os  
\!/\stabo{s}}\teta_{ur+vs},
  \qquad r,s\in \kmo.
\end{align*}
  Moreover, the generating set  $\{\teta_r\} $
is a linear basis for $\mathfrak A$, and
$\mathfrak A$ is commutative.
\item For each $d\in \ox$, the subset  $\{\teta_r:  r d \in \oo\}$  
spans a finite dimensional subalgebra
$\mathfrak A(d)$ and $\mathfrak A = \cup_{d\in \ox} \mathfrak A (d)$.

\item
The algebra $\mathfrak A$ embeds in the Hecke C*-algebra, where it
generates a C*-subalgebra $A_\theta$ that is
universal  for the above relations (interpreted now in the category of  
C*-algebras).
\item The map
\[
  \teta_r \mapsto \frac{1}{R(r)} \sum_{u \in \os\!\!/\stabo{r}}  
\delta_{ur} \ \in C^*(\kmo).
\]
determines an embedding of  $\mathfrak A$ as a *-subalgebra of $\C  
[\kmo]$ and
of $A_\theta$ as a C*-subalgebra of  $C^*(\kmo)$.
\end{enumerate}
\end{proposition}
\begin{proof}
The first two relations are immediate to verify.
In order to prove the third one,
let $\stabo{r,s}$ be the subgroup of units that fix both
$r$ and $s$ modulo $\oo$.
Clearly $\stabo{r,s} =\stabo{r\oo +s\oo}$, so it is of finite index in  
$\os$ and by \lemref{L: er is sum},
\begin{gather}\label{E:furevs}
\teta_r\convolve\teta_s=\mathfrac{1}{[\os  
:\stabo{r,s}]}\mathfrac{1}{[\os :\stabo{r,s}]}
\sum_{u\in\os \!/\stabo{r,s}}\sum_{v\in\os  
\!/\stabo{r,s}}\ER_{ur}\convolve \E_{vs}.
\end{gather}
Using \lemref{L: convolution of e and f},
we write the indicator function $\ER_{ur}\convolve \E_{vs}$ as a sum of  
indicator functions,
\[
\ER_{ur}\convolve \E_{vs}
=\sum_{w\in\os}\left[\left(\begin{matrix}1&urw+vs+\oo\\  
0&w\end{matrix}\right)\right] .
\]
In the triple summation resulting in \eqref{E:furevs}, we first sum  
over $u$ and $v$
replacing $u$ by $u/w$ when summing over $u$, to obtain
\[
\teta_r\convolve\teta_s = \mathfrac{1}{[\os :\stabo{r,s}]^2}
\sum_{w\in\os}\sum_{u,v\in\os \!/\stabo{r,s}}
\left[\left(\begin{matrix}1&\oo\\ 0&w\end{matrix}\right)
\left(\begin{matrix}1&ur+vs\\ 0&1\end{matrix}\right)\right].
\]
Fixing now the values of $r,s,u,v$, the sum over $w$ gives the  
indicator function $\E_{ur+vs}$, so that
$$
\teta_r\convolve\teta_s = \mathfrac{1}{[\os  
:\stabo{r,s}]^2}\sum_{u,v\in\os \!/\stabo{r,s}}\E_{ur+vs}
=\mathfrac{1}{[\os :\stabo{r,s}]^3}\sum_{w\in\os  
\!/\stabo{r,s}}\sum_{u,v\in\os \!/\stabo{r,s}}\E_{ur+vs},
$$
where we have introduced an extra sum over $w$.
Now, in the sum over $u$ and $v$, we replace $u$ by $wu$ and $v$ by  
$wv$ to obtain
\[
\teta_r\convolve\teta_s =\mathfrac{1}{[\os :\stabo{r,s}]^3}
\sum_{w\in\os \!/\stabo{r,s}}\sum_{u,v\in\os  
\!/\stabo{r,s}}\E_{w(ur+vs)}.
\]
Finally, we apply \lemref{L: er is sum} to the sum over $w$ to obtain
\[
\teta_r\convolve\teta_s =\mathfrac{1}{[\os  
:\stabo{r,s}]}\mathfrac{1}{[\os :\stabo{r,s}]}
\sum_{u,v\in\os \!/\stabo{r,s}}\teta_{ur+vs}.
\]
It is clear from this that the $\teta_r$ commute with each other, and   
that their
linear span is multiplicatively closed.
  We have already seen that $\teta_r = \teta_s$ if and only if $r$ and  
$s$
are in the same $\os$-orbit in $\kmo$, so as a set, $\{\teta_r : r\in  
(\kmo)^{\os}\}$
  is  linearly independent, because the supports of any two distinct  
elements are disjoint.
It follows that this set is a linear basis for $\mathfrak A$.
Conversely, if $r$ and $s$
are in the same $\os$-orbit in $\kmo$, then  the second relation  
implies that
the corresponding universal generators must coincide as well. Thus
the canonical homomorphism
of the universal unital *-algebra with the given presentation  onto  
$\mathfrak A$
maps a spanning set one-to-one and onto a linear basis; this implies  
that
the map is an isomorphism and finishes the proof of part (1).

It follows easily from the relations that $\mathfrak A(d) := \lspa\{  
\teta_r:  dr \in \oo\}$
  is a unital *-subalgebra of $\mathfrak A$, of dimension at most
$| (1/d)\oo/ \oo| = | \oo / (d\oo) | = N_d$.
Given a finite set $F = \{r_1,r_2, \ldots,r_n \} \subset \kk$,
  we choose  $d \in \ox$ such that
  $dr_i \in \oo$  for each $i = 1, 2, \ldots, n$;  then the  
$\teta_{r_i}$'s are contained in
  $\mathfrak A(d)$. This proves part (2).

To deal with the corresponding question at the C*-algebra level observe  
first
that each $\teta_r$ is self-adjoint and that it has finite spectrum  
(because it generates a finite dimensional *-subalgebra). Thus,  every  
$\theta_r$ has uniformly bounded norm
in any representation, and it makes sense to consider the universal
C*-algebra $ A^u_\theta$ of the relations.
Suppose $r$ and $s$ determine different $\os$-orbits in $\kmo$,
 i.e.\ $(r+\oo)\os \neq (s+\oo)\os$.
Using the left regular representation $\lambda$ of $\hcal(\pk,\po)$
on $\ell^2(\po\backslash \pk)$, we obtain operators $\lambda(\teta_r)$
  and $\lambda(\teta_s)$. When these operators
  act on the vector $[\po] \in \ell^2(\po\backslash \pk)$
they produce vectors with disjoint supports.
  Hence, the collection of operators $\lambda(\teta_r)$,
  indexed by (representatives of classes in) $(\kmo)^{\os}$, is linearly  
independent.
It follows that each of the finite dimensional
subalgebras $\mathfrak A(d)$ is represented faithfully in  
$\overline{\lambda(A_\theta)}$ and
  hence the canonical homomorphism
of the universal C*-algebra $A^u_\theta$ onto  
$\overline{\lambda(A_\theta)}$ is an isomorphism.
  Notice in passing that  $A^u_\theta = \injlim _{d\in \ox} \mathfrak  
A(d)$ with embeddings
  given by the inclusions $\mathfrak A(d) \subset \mathfrak A(d')$ when  
$d|d'$.
  This finishes the proof of part (3).

In order to prove part (4), we write $\vartheta_r := \frac{1}{R(r)}  
\sum_{u \in \os\!\!/\stabo{r}} \delta_{ur} \in \C[\kmo]$.
  A straightforward computation in the group algebra  $\C[\kmo]$
shows that the elements $\vartheta_r$ satisfy the
relations, so the map $\teta_r \mapsto \vartheta_r$ extends to a  
homomorphism.
To prove that this homomorphism is in fact injective, we need to show  
that the
  $\vartheta_r$'s are linearly independent.
  Suppose that $ \sum_r \lambda_r \vartheta_r$ is a vanishing finite  
linear combination,
  which we assume has been reduced so that different $r$'s come from  
different orbits in $(\kmo)^{\os}$, and so
  write as $\sum_r \sum_u\frac{\lambda_r}{R(r)} \delta_{ur}$.
  As a result of the reduction, the $\delta_{ur}$ are different, hence  
linearly independent, and
  all the $\lambda_r$ must vanish. This proves that $\mathfrak A$ embeds  
into $\C[\kmo]$.

Finally, since the finite dimensional subalgebras $\mathfrak A(d)$ are  
mapped injectively, the C*-algebra homomorphism of $A^u_\theta$ into  
$C^*(\kmo)$ determined by $\teta_r \mapsto \vartheta_r$  is injective  
too.
\end{proof}

\begin{proposition}\label{P: hecke-algebra-presentation}
The Hecke algebra $\hcal(\pk,\po)$ is (canonically isomorphic to)
the universal unital *-algebra over\/ $\C$ generated by
elements $\{\mu_a : a \in \ox \}$ and\/ $\{\teta_r : r\in \kmo\}$
subject to the relations:
\begin{enumerate}
\smallskip
\item[(I.1)] \ \  $\mu_w=1$, for all $w\in \os$,

\smallskip
\item[(I.2)] \ \ $\mu_a^*\mu_a=1$, for all $a\in \ox$,

\smallskip
\item[(I.3)] \ \ $\mu_a\mu_b=\mu_{ab}$, for all $a,b \in \ox$,

\medskip
\item[(II.1)] \ \ $\teta_0 = 1$,

  \smallskip
\item[(II.2)] \ \ $\teta_{wr} = \teta_r = \teta_r^*$, for all $r \in  
\kmo$ and $w\in \os$,

\medskip
\item[(II.3)] \ \ $\displaystyle \teta_r \teta_s = \mathfrac{1}{R(r)}
\mathfrac{1}{R(s)}\sum_{u\in\os \!/\stabo{r}}\sum_{v\in\os  
\!/\stabo{s}}\teta_{ur+vs}$,
for all $r,s\in \kmo$,

\medskip
\item[(III)] \ \
$\displaystyle \mu_a \teta_r \mu_a^* = \mathfrac{1}{N_a} \sum_{b\in  
\oo\!/a\oo}
    \teta_{\mathfrac{r+b}{a}}$ , for all $a\in \ox$ and $r\in \kmo$.
\end{enumerate}
Moreover, the set $\{\mu_a^* \teta_r \mu_b : a,b \in \ox \text{ and }  
r\in \kmo\}$ is a linear basis for
$\hcal(\pk,\po)$.
\end{proposition}

\begin{proof}
The first step is to show that the elements $\mu_a$ and $\teta_r$ of  
$\hcal(\pk,\po)$ defined in
\eqref{E: mua} and \eqref{E: er} satisfy
the relations, thus giving a canonical homomorphism of the universal
algebra of the relations onto the Hecke algebra.
The relation (I.1) follows immediately from the definition,
while (II.1),  (II.2) and (II.3) are from \proref{P:Ateta  
characterization} and \remref{R: factorthru}.
Next we observe that
the double coset supporting $\mu$ is a right coset, so for $f \in  
\hcal(\pk,\po)$,
the product $N_a^{1/2} (\mu_a  \convolve f)$ is given by
\[
N_a^{1/2} (\mu_a  \convolve f)(\gamma) = f \left( \matr{0}{a\inv}  
\gamma \right) \qquad \gamma \in \pk.
\]
Once we know how to multiply any function by $\mu_a$ on the left, the
argument given in \cite[Proposition 18]{bos-con} gives
relations (I.2) and (I.3), and since taking adjoints gives  
multiplication by
$\mu_a^*$ on the right, a further straightforward computation gives  
(III).

To avoid confusion, we shall denote by
$\tilde\mu_a$ and $\tilde\teta_r$ the universal generators subject to  
the given relations.
The rest of the proof  depends on the following key consequence of the  
relations:
\begin{gather}\label{E: extra relation}
\tilde\teta_r \tilde\mu_a = \tilde\mu_a \tilde\teta_{ar} \qquad a\in  
\ox, \ r\in \kmo.
\end{gather}
We will prove $ \tilde\teta_r \tilde\mu_a \tilde\mu_a^*  = \tilde\mu_a  
\tilde\teta_{ar} \tilde\mu_a^*$
instead; that it is equivalent to \eqref{E: extra relation} can easily be
verified on multiplying both sides on the right by the isometry $\tilde\mu_a$
or its adjoint.
We first use relations (III) and (II.3) to get
\begin{align*}
\tilde\teta_r \tilde\mu_a \tilde\mu_a^*
  & =  \tilde\teta_r \big( \mathfrac{1}{N_a} \sum_b  
\tilde\teta_{b/a}\big)  \\
& = \mathfrac{1}{N_a} \sum_b \mathfrac{1}{R(r)} \mathfrac{1}{R(b/a)}
\sum_{u\in\os \!/\stabo{r}}\sum_{v\in\os \!/\stabo{b/a}}  
\tilde\teta_{ur+ v\mathfrac{b}{a}}\\
& = \mathfrac{1}{N_a} \sum_b \mathfrac{1}{R(r)} \mathfrac{1}{R(b/a)}
\sum_{u\in\os \!/\stabo{r}}\sum_{v\in\os  
\!/\stabo{b/a}}\tilde\teta_{v(v\inv ur + \mathfrac{b}{a})}.
\end{align*}
For each fixed $v\in \os\!/\stabo{b/a}$ the new variable $w := v\inv u$  
runs over the classes of
$\os\!/\stabo{r}$, so we add first on $w$ and then simplify the  
$v$--average, because the resulting sum
does not depend on $v$:
\begin{align*}
\tilde\teta_r \tilde\mu_a \tilde\mu_a^*
& =  \mathfrac{1}{N_a} \sum_b \mathfrac{1}{R(r)} \mathfrac{1}{R(b/a)}
\sum_{v\in\os \!/\stabo{b/a}}\sum_{w\in\os \!/\stabo{r}}  
\tilde\teta_{(wr + \mathfrac{b}{a})}\\
& = \mathfrac{1}{N_a} \sum_b \mathfrac{1}{R(r)}
\sum_{w\in\os \!/\stabo{r}} \tilde\teta_{(wr + \mathfrac{b}{a})}\\
& = \mathfrac{1}{N_a} \sum_b
\mathfrac{1}{R(r)}\sum_{w\in\os \!/\stabo{r}}  
\tilde\teta_{\mathfrac{war +b}{a}} \\
& = \mathfrac{1}{R(r)}\sum_{w\in\os \!/\stabo{r}} \tilde\mu_a  
\tilde\teta_{war}\tilde\mu_a^* \qquad
\text{ by (III), for each } w\\
& = \tilde\mu_a \tilde\teta_{ar}\tilde\mu_a^*
\end{align*}
because $\tilde\teta_{war} = \tilde\teta_{ar} $ for each $w\in\os$.

Next we show that the canonical homomorphism of the universal *-algebra  
with
the given presentation to $\hcal(\pk,\po)$ is an isomorphism.

Using the given relations,  \eqref{E: extra relation}, and its immediate
consequence $\tilde\mu_a ^*\tilde\teta_r \tilde\mu_a =   
\tilde\teta_{ar}$,
one proves that the linear span of the set of monomials
$\{\tilde\mu_a^* \tilde\teta_r \tilde\mu_b: a,b\in \ox , r\in \kmo\}$  
is multiplicatively closed and
$*$-closed. The computation is analogous to the one found
in the proof of \cite[Lemma 1.8]{alr}.
Since these monomials include the generators $\tilde\mu_a$ and
$\tilde\teta_r$, as well as the identity,
  they span the universal unital *-algebra with the above presentation.

At the level of the Hecke algebra, one computes that the bi-invariant  
function $\mu_a^*\teta_r\mu_b$ is
supported on the (single) double coset of $\smatr{rb}{b/a}$, where it  
takes the value
$([\stabo{ar}:\stabo{r}]\sqrt{N_{ab}})^{-1}$, 
see the proof of \cite[Theorem 2.3]{alr} for a similar computation. Since every double coset  
arises as the support of such a function for some $a,b \in \ox$
and $r\in \kmo$, the set of such products spans the  
Hecke algebra.
The {\em caveat} of \cite[Remark 1.9]{alr} applies here too, to the  
effect that there is some redundancy in
the labelling of the universal spanning monomials by the triples  
$(a,r,b)$, and,
of course, at least as much redundancy in the corresponding labelling of
the monomials spanning the Hecke algebra. We claim that this redundancy  
is exactly the same
in both cases; more explicitly, we claim
that if $\mu_a^* \teta_r \mu_b$ and $\mu_c^* \teta_s \mu_d$ do not have  
disjoint supports
as functions on $\pk$ then $\tilde\mu_a^* \tilde\teta_r \tilde\mu_b =   
\tilde\mu_c^* \tilde\teta_s \tilde\mu_d$ in the universal
algebra of the relations, and hence also $\mu_a^* \teta_r \mu_b  
=\mu_c^* \teta_s \mu_d$ in $\hcal(\pk,\po)$, so the elements $\mu_a^* \teta_r \mu_b $
are linearly independent in $\hcal(\pk,\po)$.
This will complete the proof, because it implies that the canonical  
homomorphism maps a
spanning set one-to-one and onto a linear basis, and hence is an  
isomorphism.
Recall that  the support of $\mu_a^* \teta_r \mu_b $ is the single  
double coset
$\po \smatr{rb}{b/a}\po $ and that of  $\mu_c^* \teta_s \mu_d$ is
$ \po\smatr{sd}{d/c}\po$.
If these supports are not disjoint they coincide, so
there exist $u$ and $w$ in $\os$ and  $m \in \oo$ such that
$b/a = u (d/c)$ and $rb = w (sd) + m(b/a) \bmod  \oo$.
Then, going along the lines of \cite[Remark 1.9]{alr},
\begin{align*}
\tilde\mu_a^* \tilde\teta_r \tilde\mu_b
& = \tilde\mu_a^*\tilde\mu_b \tilde\teta_{br} \qquad \text{ by  
\eqref{E: extra relation}}\\
& = \tilde\mu_a^*\tilde\mu_b \tilde\teta_{w ds + m b/a}    \\
& = \tilde\mu_c^* \tilde\mu_d \tilde\teta_{w ds + m u d/c} \qquad
\text{ because } \tilde\mu_a^* \tilde\mu_b = \tilde\mu_c^* \tilde\mu_d\\
& = \tilde\mu_c^* \tilde\teta_{ws + m u/c} \tilde\mu_d
\qquad \text{ by \eqref{E: extra relation}}\\
& = \tilde\teta_{cws + m u} \tilde\mu_c^* \tilde\mu_d
\qquad \text{ by \eqref{E: extra relation}}\\
& = \tilde\mu_c^* \tilde\teta_s \tilde\mu_d,
\end{align*}
finishing the proof of the claim and of the proposition.
\end{proof}

It will be important for our study of KMS states to have a presentation  
of the  Hecke C*-algebra as a semigroup crossed product and in terms of generators and  
relations. With this purpose, and in order to introduce the notation,
 we briefly review now the definition of a semigroup crossed
product. When $S$ is a semigroup that acts by endomorphisms $\alpha_s$ of the 
 unital C*-algebra $A$, we say that $(A, S, \alpha)$ is a semigroup dynamical system.
 A {\em covariant representation} of such a system is a pair $(\pi, v)$
 consisting of  a unital representation $\pi$ of $A$ on a Hilbert space
 and a representation $v$ of $S$ by isometries on the same Hilbert space,
 such that the covariance condition $\pi(\alpha_s(a)) = v_s \pi(a) v_s^*$ 
 is satisfied for every $a\in A$ and $s\in S$. 
 The semigroup crossed product associated to $(A,S,\alpha)$ is
 the C*-algebra generated by a universal covariant representation.
 It is unique up to canonical isomorphisms and is denoted $A\rtimes_\alpha S$.
  When the endomorphisms are injective, the component $\pi$ of
 the universal covariant representation is injective 
and it is customary to drop it from the notation and to think of $A\rtimes_\alpha S$
as being generated by a copy of $A$ and a semigroup  of isometries $\{v_s: s\in S\}$ 
 that are universal for the covariance condition.
 See \cite{quasilat} for more details on semigroup crossed products.
 The next proposition gives the appropriate semigroup  
dynamical system for our Hecke algebra.
We point out that the endomorphisms arising in this particular construction
are  {\em injective corner endomorphisms}, that is, each $\alpha_s$
 is an isomorphism of $A$ onto the corner $p_sAp_s$ 
determined by the projection $p_s = \alpha_s(1)$. See \cite{minautex} for general facts about
crossed products by semigroups of corner endomorphisms.

\begin{proposition} \label{actionalpha}
There is an action of the semigroup
$\ox/\os$ of principal integral ideals of $\oo$
by injective corner endomorphisms of the C*-algebra $A_\theta$, given by
\begin{gather}\label{E: alpha_a}
\alpha_a(\teta_r) = \mathfrac{1}{N_a}\sum_{b\in\oo\!  
/a\oo}\teta_{\mathfrac{r+b}{a}}, \qquad a\in \ox.
\end{gather}
  For each $a\in \ox$, $\alpha_a(1)=\mu_a\mu_a^*$ is a projection, and
if the ideal $a\oo \cap b\oo$ is principal, then
\begin{gather}\label{E:alfa lcm}
\alpha_a(1) \alpha_b(1) = \alpha_{a\oo \cap b\oo}(1).
\end{gather}
The action $\alpha$ has a left inverse action $\alphin$ by surjective  
endomorphisms  determined by
$\alphin_a: \teta_r \mapsto \teta_{ar}$, such that $\alphin_a\circ  
\alphin _b = \alphin _{ab}$,
$\alphin_a \circ \alpha_a = \id $ and $\alpha_a \circ \alphin_a $
is multiplication by $\alpha_a(1)$.
\end{proposition}

\begin{proof}
Notice first that the relation (II.2) makes it clear that $\alpha_a$  
depends only on the
ideal $a\os \in \ox/\os$ and not on the specific representative $a\in  
\ox$.
To show that $\alpha$ is a semigroup of endomorphisms, and to define  
$\alphin$, we use the copy of $A_\theta$ sitting inside
the universal algebra $C^*(\mu,\teta)$ of the relations, where $\alpha$  
and $\alphin$ are implemented
by the isometries:
\begin{gather}\label{alfabeta}
    \alpha_a(f) =  \mu_a f \mu_a^* \qquad \text{and}\qquad  \alphin_a(f)  
:= \mu_a^* f \mu_a \qquad \text{ for } f\in \mathcal H(\pk, \po).
\end{gather}
Since the $\mu_a$ form a semigroup of isometries, it is easy to see that
  $\alpha$ and $\alphin $ have the desired properties.
  Multiplying \eqref{E: extra relation} on the left by $\mu_a^*$
shows that $\alphin_a(\teta_r)=  \mu_a^* \teta_r \mu_a = \teta_{ar}$.

To prove \eqref{E:alfa lcm} suppose first that $a$ and $b$ are  
relatively prime;
then $a/u$ and $b/v$ are also relatively prime for every  $u,v \in \os$,
and from \cite[Proposition 1.2]{alr} we know that
\[
\big(\mathfrac{1}{N_a}\sum_{x\in\oo\! /a\oo}\delta_{\mathfrac{ux}{a}}  
\big)
\big(\mathfrac{1}{N_b}\sum_{y\in\oo\! /b\oo}\delta_{\mathfrac{vy}{b}}  
\big)
=
\mathfrac{1}{N_{ab}}\sum_{z\in\oo\! / ab \oo}\delta_{\mathfrac{z}{ab}}.
\]
Averaging first over $u\in \os\!/\stabo{1/a}$ and then over $v\in  
\os\!/\stabo{1/b}$ gives
\[
\big(\mathfrac{1}{N_a}\sum_{x\in\oo\! /a\oo}\vartheta_{\mathfrac{x}{a}}  
\big)
\big(\mathfrac{1}{N_b}\sum_{y\in\oo\! /b\oo}\vartheta_{\mathfrac{y}{b}}  
\big)
=
\mathfrac{1}{N_{ab}}\sum_{z\in\oo\! /ab\oo}\vartheta_{\mathfrac{z}{ab}},
\]
which proves
\[
\alpha_a(1) \alpha_b(1) = \alpha_{a b}(1),
\]
  in the copy of $A_\theta$ inside of $C^*(\kmo)$ given in part (4) of  
\proref{P:Ateta characterization}.
When $a$ and $b$ are not relatively prime but have a principal l.c.m.  
$a\oo \cap  b\oo = c\oo$, we write
$a' = c b\inv$ and $b' = c a\inv$, and we let $d: = abc\inv$ be the  
g.c.d. of $a$ and $b$.
Then $a'$ and $b'$ are relatively prime and
\eqref{E:alfa lcm} follows from the above because
\[
\alpha_a(1) \alpha_b(1) = \alpha_{da'}(1)  \alpha_{db'}(1) =
\alpha_d(\alpha_{a'b'}(1) ) = \alpha_c(1).
\]
\end{proof}

\begin{corollary}
If\/ $a\oo + b\oo$ is principal, so that $a\oo + b\oo=c\oo$ for some   
$c\in \ox$,
and if we write $a = a'c$ and $b = b'c$ for $a', b' \in \ox$, then
\begin{gather}\label{E: extracommrelation}
\mu_b^* \mu_a = \mu_{a'} \mu_{b'}^*.
\end{gather}
\end{corollary}
\begin{proof}
Multiply \eqref{E:alfa lcm} on the left by $\mu_a^*$ and on the right  
by $\mu_b$, and then use \eqref{alfabeta} with $f = 1$ and the fact that 
$a\oo \cap b\oo = a' b' c \oo$.
\end{proof}

\begin{theorem}\label{T: universal presentation}
The Hecke C*-algebra $\heck$ is canonically isomorphic to the semigroup  
crossed product
$A_\theta \rtimes_\alpha (\ox\!/\os)$
and to the universal unital C*-algebra $C^*(\tilde\mu, \tilde\teta)$  
with presentation
(I), (II), (III) from \proref{P: hecke-algebra-presentation}.

Similarly, the Hecke algebra is canonically isomorphic to the  
{`algebraic'} semigroup crossed product $\mathfrak A \rtimes_\alpha  
(\ox\!/\os)$, which embeds as the *-subalgebra of
$A_\theta \rtimes_\alpha (\ox\!/\os)$ spanned by the
monomials $\mu_a^* \teta_r \mu_b$.
\end{theorem}

\begin{proof}
The isomorphism between the semigroup crossed product and the  
C*-algebra of the given presentation is an easy consequence of  
\proref{P:Ateta characterization} and the
universal property of semigroup crossed products \cite{quasilat}.

We have already seen that the operators $\lambda(\teta_r)$
generate a faithful representation of $A_\theta$ on  
$\ell^2(\po\backslash \pk)$ and that the
$\lambda(\mu_a)$ form a semigroup of isometries.
Moreover, the relation (III) says that the pair $(\lambda|_{A_\theta},  
\lambda\circ \mu)$  is a covariant representation
of the semigroup dynamical system $(A_\theta, (\ox\!/\os), \alpha)$,
  so the universal property of $A_\theta \rtimes_\alpha (\ox\!/\os)$
gives a homomorphism of $A_\theta \rtimes_\alpha (\ox\!/\os)$ onto
$\heck$, which extends the representation $\lambda$ of $\hcal(\pk,\po)$
as convolution operators on the space $\ell^2(\po\backslash \pk)$. We  
will
refer to this as the {\em Hecke representation} of the Hecke
  algebra and C*-algebra,  and denote it by $\lambda$;
our next goal is to show that it is injective.
To do this, we will show that the vector state $\omega_{[\po]} :=  
\langle \lambda(\cdot)
[\po] , [\po]\rangle$  corresponding  to the vector $[\po] \in  
\ell^2(\po\backslash \pk)$ is a faithful state
of the crossed product. From this it will follow that $\lambda$ is  
faithful
because its cyclic subrepresentation associated to $[\po]$ is faithful.

In order to see that $\omega_{[\po]} $ is faithful,
  let $E_{\hat{\alpha}} : A_\theta \rtimes_\alpha (\ox\!/\os) \to   
A_\theta$ be
the conditional expectation of the dual action of $\widehat{\kk^*}$,  
obtained by
  averaging over the compact $\widehat{\kk^*}$--orbits.
Using the embedding of $A_\theta$ in $C^*(\kmo)$ from part (4) of
  \proref{P:Ateta characterization},
  restrict the canonical trace on $C^*(\kmo)$ to a trace $\tau$ on  
$A_\theta$ .
Since $\tau$ and $E_{\hat{\alpha}}$ are faithful positive maps (in the  
sense that
their kernels contain no nontrivial positive element), their composition
$\tau \circ E_{\hat{\alpha}}$ is a faithful state on $A_\theta  
\rtimes_\alpha (\ox\!/\os)$.
Next observe that  the state $\omega_{[\po]}$
factors through $E_{\hat{\alpha}}$ and coincides with $\tau$ on the  
$\teta_r$'s.
Since these generate the range of
$E_{\hat{\alpha}}$, we have $\tau \circ E_{\hat{\alpha}} =  
\omega_{[\po]} $, finishing the proof.
See \cite[Example 1.12 and proof of Proposition 2.4]{alr} for a similar  
computation carried out in detail.
\end{proof}

It is possible to give a different proof of the semigroup crossed  
product structure of $\heck$
using results from \cite{LaLa,LaLaErr}. This alternative
approach reveals that $A_\theta$ is itself a Hecke C*-algebra,
although it does not yield the explicit formula for products (II.3)  
obtained above by direct methods, so it paints a complementary picture.
In order to apply \cite[Theorem 1.9]{LaLa} we first need to write $\pk$  
as a
convenient semidirect product,
for which we need a multiplicative cross section of the quotient map
  $\ks \longrightarrow \ks /  \os$ that maps $\ox/\os$ back into $\ox$.
  When every ideal in $\oo$ is principal, such a cross section is easy  
to obtain; one
simply selects a prime generator for each prime ideal and then extends
the map freely and multiplicatively.
The existence of such cross sections
when the ideal class group is nontrivial is a bit more delicate.

\begin{lemma}\label{L: section}
The quotient map $\ks \longrightarrow \ks /\os $ has a multiplicative  
cross section
mapping $\ox/\os$ into $\ox$.
\end{lemma}

\begin{proof}
Write the ideal class group of $\kk$ as a product of cyclic groups
  $C_{d_1}\times\dots\times C_{d_r}$,
with $d_1\divides\dots\divides d_r$.
Choose prime ideals $\pid_1,\dots ,\pid_r$ mapping to the generators of
  $C_{d_1}, \dots,C_{d_r} $.
Thus a product $\prod_j\pid_j^{e_j}$ is a principal ideal if and only  
if $d_j\divides e_j$ for each $j = 1, \dots, r$.
Choose elements $a_j\in\oo$ such that $\pid_j^{d_j}=a_j\oo$.
Also choose, for each prime ideal $\pid$, (including the $\pid_j$'s) a  
product
$\prod_j\pid_j^{e_j(\pid )}$ in the same ideal
class, and finally, choose an element $a_\pid\in \kk$ such that
  $\pid =a_\pid\prod_j\pid_j^{e_j(\pid )}$.

Let $\aid$ be a principal ideal of $\oo$, and factor $\aid$ as a  
product $\prod_\pid\pid^{e_\pid(\aid)}$ of prime ideals. Replacing the  
above expression for $\pid$
  gives
   \[
  \textstyle
  \aid = \prod_\pid  
a_\pid^{e_\pid(\aid)}\cdot\prod_\pid\prod_j\pid_j^{e_\pid(\aid)  
e_j(\pid )}
  =\prod_\pid a_\pid^{e_\pid(\aid)}\cdot \prod_j\pid_j^{\sum_\pid  
e_\pid(\aid) e_j(\pid )}.
\]
Since $\aid$ is principal,
it follows that $\sum_\pid e_\pid(\aid) e_j(\pid )$ is a multiple of  
$d_j$, for each $j= 1, \dots , r$,
and we obtain a generator of the ideal $\aid$
\[\textstyle
\aid =\prod_\pid a_\pid^{e_\pid(\aid)}\cdot\prod_j  
a_j^{(1/d_j)\sum_\pid e_\pid(\aid) e_j(\pid )}\oo .
\]
We choose this generator of $\aid$ as the image of $\aid\in\ox /\os$ in  
$\ox$.
First of all, since $\aid$ is an ideal of $\oo$, this generator is an  
integer.
And secondly, since $e_\pid(\aid \bid) = e_\pid(\aid) + e_\pid( \bid) $
  for every prime ideal $\pid$, and fractional ideals $\aid$ and $\bid$,
  it is easily verified that this construction gives a multiplicative  
section
  $\ks/\os \to \ks$ that maps
$\ox/\os$ into $\ox$.
\end{proof}

  \begin{corollary}
The subgroup $N := \smatr{\kk}{\os}$ is normal in  $\pk $ and
contains $\po$ as a Hecke subgroup.
The Hecke *-algebra  of the `intermediate' Hecke inclusion $\po \subset  
N$
is canonically isomorphic to the *-algebra $\mathfrak A$ with  
presentation given in
\proref{P:Ateta characterization}. Along the same lines, the Hecke  
C*-algebra
  of this inclusion is the C*-algebra $A_\theta$ with the same  
presentation.
Moreover, we have canonical isomorphisms at the level of *-algebras and  
C*-algebras:
\begin{align*}
\mathcal H (\pk,\po) &\cong \mathcal H(N,\po) \rtimes(\ox/\os) \quad  
\text and
\\
\heck &\cong C_r^*(N,\po) \rtimes (\ox/\os).
\end{align*}
\end{corollary}
\begin{proof}
Since $\po$ is Hecke in $\pk$, it is obviously Hecke in $N \subset  
\pk$, and
$N$ is normal because it is the kernel of
the homomorphism
$\smatr{y}{x} \in \pk  \mapsto x \os \in \ks /\os$.
  Let $\sigma $ be a cross section given by the preceding lemma
and consider the map $\tilde\sigma: \ks/\os \to \pk$, defined by   
$\tilde\sigma(x\os) =  \smatr{0}{\sigma(x\os)}$. Using this cross  
section, we may view $\pk$ as the semidirect product
$N \rtimes (\ks/\os)$, and \cite[Theorem 1.9]{LaLa} does the rest.
\end{proof}
\begin{remark}1) Although the cross section $\sigma$ is noncanonical,  
the resulting endomorphisms of the small Hecke algebra $\hcal(N,\po)$  
are independent from it,
essentially because the double cosets in $N$ are $\os$--invariant.

2) Because Lemma 1.3 of \cite{LaLa} is wrong as stated, 
a correction factor is necessary on the spanning monomials of
  \cite[Theorem 1.9 (ii)]{LaLa} for the set to be a linear basis,
 see~\cite{LaLaErr}.  In the present case, the factor is $1/R(r)$,
  and is already naturally included in the definition of the monomials
  in \eqref{E: er}.
A very interesting generalization of \cite{LaLa}
to nonsplit extensions, as well as the correct substitute to the  
incorrect lemma, can be found in \cite{BigDog}.
\end{remark}

\section{The adelic semigroup crossed product}
\label{S: alternative}
It is very convenient to
have a realization of the Hecke algebra $\heck$ as a semigroup crossed  
product
on which it is possible to define the endomorphisms and the group of  
`geometric' symmetries
as transformations of a compact space. We shall obtain this
realization via the Fourier--Gelfand transform
of $\kmo$, for which we need to review first the self-duality pairing of  
the additive group of full adeles.
   We refer to \cite{Tate} or \cite{lan} for the details.
Denote by $\mk$ the set of equivalence classes of valuations on
$\kk$, with $\mko$ the set of finite (nonarchimedean) valuations, and  
let $\kk_v$
and $\oo_v$ denote the completion of $\kk$ and $\oo$ with respect to  
the valuation $v$.
The ring  of full adeles over $\kk$ is, by definition, the restricted  
product
$\A (\kk):= \prod_{v\in\mk} (\kk_v; \oo_v)$; the additive group is self  
dual. We shall consider the pairing of
$\A (\kk)$ with itself given in \cite{Tate}:
  suppose first $p\in \N$ is a prime, and let $\chi_p$ be the character  
of $\Q_p$ given by
\[
\chi_p(x)=\exp (2\pi i \lambda_p (x)), \qquad x\in \Q_p,
\]
where $\lambda_p(x)\in\Q$ is chosen so that  $x-\lambda_p(x)\in\Z_p$  
and is well-defined modulo $\Z$.
Also for $p = \infty$ define a character $\chi_\infty$ of $\Q_\infty =\R$ by
\[
\chi_\infty (x)=\exp (-2\pi ix), \qquad x\in \Q_\infty.
\]
For each adele $x=(x_\infty ,x_2,x_3,x_5, \dots ) \in \A(\Q)$, define
\[
\chi_{\A (\Q )}(x)=\prod_p\chi_p(x_p).
\]
Clearly, $\chi_{\A(\Q)}(x)=1$ for $x\in\Q$.
  Let $\Tr_{\A(\kk)/\A(\Q)}(a)=\bigl(\sum_{v|p}\Tr_{\kk_v/\Q_p}(a_v)\bigr)_{p}$
  denote the trace of an adele $a$ over $\kk$ to the adeles over $\Q$,
where $p=\infty,2,3,\dots$ runs over the valuations of $\Q$.
For adeles $a,b\in\A (\kk )$,
we define the canonical pairing
 $\left< a,b\right> =\chi_{\A (\Q )}(\Tr_{\A (\kk )/\A (\Q )}(ab))$, and
  this gives an isomorphism of $\A(\kk)$ onto its Pontrjagin dual  
$\widehat{\A(\kk)}$,
  in which the adele $b$ is mapped to the character on $\A(\kk)$
  given by $a \mapsto   \left< a ,b\right> $. Once this duality has been  
established,
  general considerations give pairings for subgroups and quotients of  
$\A(\kk)$.

For instance, consider $\kk \subset \A(\kk)$,
embedded diagonally,
 and let~\mbox{$\kk^\perp$} be the
 group~\mbox{$\kk^\perp :=\{ b\in \A(\kk): \left< a,b\right> =1 \text{ (all }a\in\kk)\}$}
of characters that are trivial on $\kk$.
It is easy to see that $\kk\subset \kk^\perp$, and the reverse inclusion
holds by~\cite[Theorem~4.1.4]{Tate}.
Thus the duality pairing of $ \A (\kk )$ with itself  determines an  
isomorphism of the dual of $\kk$ to the quotient $ \A (\kk )/\kk$.
This, in turn, yields an isomorphism of the dual of the quotient $\kk  
/\oo$ to
the subgroup  of $\A (\kk )/\kk$
consisting of those adeles (modulo $\kk$) that induce the trivial  
character on $\oo$. Specifically,
define
\[
\oo^\perp := \{ b \in \A (\kk ): \left<a,b\right> = 1 \text{ for all }  
a\in \oo\} ;
\]
then there is a  duality pairing of the classes $k +\oo \in \kmo$
and $b+\kk \in \oo^\perp/\kk$ given by $\left<k+\oo , b+\kk\right> :=  
\left<k, b\right>$.

This is not quite enough to see
endomorphisms and symmetries
as coming from multiplication; we need a description
of $\oo^\perp/\kk$ in terms of a subset of the finite adeles,
on which multiplication is defined.
Suppose first $v$ is a nonarchimedean valuation on $\kk$,
and let $p$ be the rational prime such that $v|p$.
The local version of the duality, pairing $\kk_v$ with itself, is
given by  $\left< a,b\right>_v := \chi_p(\Tr_{\kk_v/\Q_p}(ab))$ for  
$a,b \in \kk_v$.
Since $\left< a,b\right>_v = 1$ for all $a\in \oo_v$ if and only if
$\Tr_{\kk_v/\Q_p}(ab) \in \Zp$ for all $a\in \kk_v$, the set of
elements that induce the trivial character on $\oo_v$ is
  \[
  \invdv = \{ b\in\kk_v\colon \Tr_{\kk_v/\Q_p}(ab) \in \Z_p \text{ for  
all } a\in \oo_v\},
  \]
which is classically referred to as the {\em local inverse different}  
at $v$.
The set $ \invdv $ is a fractional ideal in $\kk_v$
that clearly contains $\oo_v$, so that its inverse, the {\em local  
different} $\diff_v$
is an integral ideal in $\kk_v$.
If we denote by  $\mathfrak P_v$ the maximal ideal of $\oo_v$, then  
$\diff_v = \mathfrak P_v^{d_v}$
for some $d_v \geq 0$. One proves that the support of $d$ is finite  
using the approximation theorem;
  in fact,  $\diff_v = \oo_v$ if and only if $v\in \mko$ is unramified.

The space we need is the product of the local inverse differents:
\[
\invd := \prod_{v\in \mko}  \invdv = \{ b\in \A^0(\kk) : \chi_{\A (\Q  
)}(\Tr_{\A (\kk )/\A (\Q )}(ab)) = 1 \ \forall a\in
{\textstyle\prod_{v\in \mko}}  \oo_v\},
\]
which is an additive subgroup  of the finite adeles $ \A^0(\kk)$, but
which we view as a subset of the full adeles by attaching a trivial  
component
at every infinite place. By abuse we will refer to $\invd$ as the  
(adelic) inverse different,
observing that it is a cartesian product $\invd = \prod_v \mathfrak  
P_v^{-d_v}$.

In order to avoid confusion, we point out that in the standard  
terminology the
{\em global inverse different} is the set
$$
\diff\inv:= \{x\in \kk : x\in \invdv \text{ for all finite } v\}.
$$
It is a fractional ideal in $\kk$, whose prime factorization is
$\diff\inv = \prod_{\mathfrak P} \mathfrak P^{-d_v}$, where
$d_v$ is as above for the valuation corresponding to the prime ideal 
$\mathfrak P = \mathfrak P_v \cap \oo$.

\begin{proposition}\label{dualitypairing}
The canonical duality pairing of $ \A (\kk )$ to itself
induces an isomorphism of the dual of $\kk /\oo$ to
$\invd \times \prod_{v|\infty} \{0\}$.
\end{proposition}
\begin{proof}
 From the definition of $\invdv$ it is clear that the elements of  
$\prod_{v\in \mko}  \invdv $,
viewed as full adeles with trivial components at the
infinite places,  induce the trivial character on $\oo$ via the duality  
pairing.
Thus, we have that $\prod_{v\in \mko}  \invdv \times \prod_{v|\infty}  
\{0\} \subset \oo^\perp$.
We need to show that this inclusion induces an isomorphism of
$\prod_{v\in \mko}  \invdv $ to $\oo^\perp/\kk$.

By the Approximation Theorem,
every adele can be written, modulo $\kk$,  as $d+x$,
where $d\in\prod_{v\in \mko}  \invdv \times  \prod_{v|\infty} \{0\} $  
and
$x\in \prod_{v\in \mko} \{0\} \times \prod_{v\divides \infty}\kk_v$.
Then
\[
\A (\kk )/\kk \cong\prod_{v\in \mko}  \invdv \times \big(  
{\textstyle\prod}_{v\divides \infty}\kk_v \big) / \diff\inv.
\]
Since every nontrivial element of
$\textstyle \prod_{v\divides \infty}\kk_v/\diff\inv$ gives a nontrivial character on
$\oo$, it follows that $\prod_{v\in \mko}  \invdv \cong \oo^\perp/\kk$.
\end{proof}

\begin{definition}\label{D: extreme point}
An {\em extreme point\/} of the inverse different (or an {\em extreme  
inverse different})
  is an element $\chi$ of $ \invd  $ such that
  $\chi\rr = \invd $; equivalently, $\chi$ has  maximal valuation  
(minimal exponents) within
  $\prod_{v\in \mko}  \invdv $ at every place.
   The set of extreme points will be denoted by $\partial \invd$.
  \end{definition}
The extreme inverse different corresponds to the
subset $\mathcal X_\kk$ of $\widehat\kmo$ constructed in  
\cite[Corollary 3.5]{alr}
using elementary Fourier analysis. The extreme points play the role of  
the appropriate power of
a uniformizing element in a local field extension.
For each extreme inverse different $\chi $ the map $ w \mapsto \chi w $
gives a one-to-one correspondence:
\[
\rr := \prod_{v\in \mko} \ov \quad \longrightarrow \quad \invd  
:=\prod_{v\in \mko}  \invdv,
\]
such that
\[
W:= \prod_{v\in \mko} \ov^* \quad \longrightarrow \quad \partial \invd  
= \partial (\prod_{v\in \mko}  \invdv ).
\]

We may now write the Gelfand--Fourier transform of $C^*(\kmo) $
as the isomorphism $C^*(\kmo) \cong C(\invd)$ determined by
the pairing of \proref{dualitypairing}, namely  $\delta_r\mapsto  
\hat\delta_r = \left< r, \cdot\right>$
for $r\in \kmo$.
We note in passing that the various duality pairings of $\kmo$ to  
$\mathcal R$ used in
\cite{har-lei,alr,coh} involve a
noncanonical choice of inverse different $\chi$ (or of a point in  
$\mathcal X$
in the case of \cite{alr}). In our notation, the transforms from  
\cite{har-lei,alr,coh}
  would be given by $\delta_r \mapsto \left<  r, \chi\, \cdot \,\right>$.

To describe the symmetries of our dynamical system,
we need to analyze the action of the group $W:= \prod_{v\in \mko} \ov^*$
on the inverse different, and in particular the action of the subgroup  
$\os$ of units,
  diagonally embedded in $W$.
By Dirichlet's unit theorem, $\os$  is infinite unless $\kk$ is  
quadratic imaginary,
but since $W$ is compact,  the closure $\osbar$ of the diagonal
copy of $\os$ in $W$ is a compactification of $\os$. Next we see that  
this compactification
coincides with the natural profinite limit determined by the action of  
$\os$ on $\kmo$,
and, by duality, with the transposed action on $\invd$.
We recall that  $\stabo{r}$ denotes the stability
subgroup of $r\in \kmo$ for the action of $\os$ by multiplication on  
$\kmo$.

\begin{lemma}\label{aside}
For each $a\in \ox$ let $G_a := \os/\stabo{1/a}$. The collection  
$\{G_a\colon a\in \ox\}$,
endowed with the natural homomorphisms
\[
u \stabo{1/b} \in G_b  \mapsto u \stabo{1/a} \in G_a \qquad \text{for  
$u\in \os$, and $a|b$ in $\ox$}
\] is a projective system of finite groups.
The quotient maps $u \mapsto u \stabo{1/a}$ of\/ $\os \to G_a$  
determine an
embedding of\/~$\os$ as a dense subgroup of the  profinite group  
$\projlim G_a$,
which extends to the closure $\osbar$ of (the diagonally embedded copy  
of)\/~$\os$
inside $W$, giving an isomorphism $\osbar \cong \projlim G_a$.
\end{lemma}

\begin{proof}
The collection of subgroups $ \stabo{1/a}$ for $a\in \ox$ form an  
injective system
  (with~$\ox$ directed by divisibility) and hence the corresponding  
quotients $G_a$ form a projective
system.

Next observe that for any two distinct units $u$ and $v$ in $\os$,  
there is
$b \in \ox$ such that $u-v \notin b\oo$, hence $u/b \neq v/b \bmod \oo$,
proving that  the canonical map of $\os$ to $\projlim G_a$ is injective.

The group of units $\os$ maps onto each of the $G_a$'s by definition;
  we need to extend this map to the closure of $\os$ in $W = \prod_{v\in  
\mko}\oo_v^*$, for which
we will use that $W$ itself can be viewed as the limit
of the projective system $\{(\oo/a\oo)^* \colon a \in\ox\}$, where  
$(\oo/a\oo)^*$ is the unit
  group of the finite ring $(\oo/a\oo)$.
Suppose now $u_\lambda \to u$ in $W$. Then $(u_\lambda u\inv) \to 1$
in $(\oo/a\oo)^*$, so eventually $(u_\lambda u\inv) = (u_\mu u\inv)$
for some fixed $\mu$.
It follows that  $u_\lambda = u_\mu \bmod a\oo$, that is,
$(u_\lambda - u_\mu ) / a \in \oo$. Hence $u_\lambda u_\mu\inv$ fixes  
$1/a$ modulo $\oo$,
so $u_\lambda$ and $ u_\mu$ determine the same element of $G_a$. This  
implies that the
quotient map $\os \to (\oo/a\oo)^*$ extends to a surjective  
homomorphism $h_a: \osbar \to (\oo/a\oo)^*$ for each $a\in \ox$.

It is now easy to verify that these homomorphisms $h_a$ form
a coherent family of continuous surjections of $\osbar$
to the projective system of groups $G_a$. The intersection of
the kernels is trivial, and by the universal property of the projective  
limit
there is an embedding of $\osbar$ in $\projlim G_a$. Since the range is  
dense
and $\osbar$ is compact, this embeding is surjective and we have the  
isomorphism
  $\osbar \cong \projlim G_a$.
\end{proof}

It is easy to see that multiplication by integral ideles $W$ is  
continuous and
leaves the inverse different invariant,
thus the diagonal embedding $\os \hookrightarrow W$ determines a  
multiplicative
  action of $\os$ on $\invd$.
At the level of the C*-algebra $C(\invd)$, a function is fixed by the  
action of $\os$ if and only if
it is fixed by $\osbar$.
Averaging over the $\osbar$--orbits using the normalized Haar measure  
gives a faithful
positive conditional expectation $E_{\osbar}: C(\invd) \to  
C(\invd)^{\osbar}$.
We shall denote by $\Omega$ the (quotient) space of $\osbar$--orbits in  
$\invd$; it is a
compact Hausdorff space such that $ C(\invd)^{\osbar} \cong C(\Omega)$.
Notice that $\os$ acts trivially on $\Omega$, so the multiplicative action of $W$ on
$\invd$ drops to an action of $W/\osbar$ on $\Omega$,
which we also write as multiplication. Similarly,
 multiplication of an orbit $x\osbar$ in $\Omega$
by a principal ideal $a\in \ox/\os$ is well defined, since $\os$ has  
been factored out from everything in sight.

\begin{proposition}\label{dual+faith}
Let $\Omega$ be the orbit space of the action of $\osbar$ on $\invd :=  
\prod_v\invdv$.
There is an action $\alpha$ of\/ $\ox/\os$ by injective endomorphisms  
of\/ $C(\Omega)$ defined,
for $a\os \in \ox/\os$ and $f\in C(\Omega)$,  by
\begin{gather}\label{E: alter-def}
\alpha_a(f)(x)  : =   \left\{
\begin{array}{ll}  f(a\inv x) &\text{ if } x \in a  \Omega\\
                   0 & \text{ if } x \notin a  \Omega.
\end{array}
\right.
\end{gather}
A representation of the crossed product $C(\Omega) \rtimes_\alpha  
(\ox/\os)$ is faithful if and only if it is faithful on $C(\Omega)$.
\end{proposition}

\begin{proof}
One first defines $\alpha_a(f)$ for $f\in C(\invd)$ and $a\in \ox$ and  
then restricts to
$C(\Omega) \cong C(\invd)^{\osbar}$, where units act trivially, to  
obtain the action of  $\ox/\os$.
This shows that \eqref{E: alter-def} is independent of the point $x$  
representing an orbit and
of the integer $a$ representing a principal ideal.
Recall from the discussion preceding \proref{actionalpha}
 that the semigroup  crossed product $C(\Omega) \rtimes_\alpha  
(\ox/\os)$ is generated by a copy of $C(\Omega)$, and a
semigroup of isometries $\{v_a: a\in \ox/\os\}$.
To prove the faithfulness statement,
 let $F$ be a finite subset of $\ox/\os$ and consider the linear combination
  $\sum_{a,b\in F} f_{a,b} v_a^* v_b$;
  such elements span the associated `algebraic crossed product', which is a  
dense $*$-subalgebra of $C(\Omega) \rtimes
(\ox/\os)$.  Choose representatives in $\ox$ for the elements of $F$,  
and let
$q_1 \in C(\invd)$ be the projection constructed as in the proof of  
\cite[Lemma 4.3]{alr}. Since this projection
  is invariant under $\os$, it lies in $C(\Omega)$. The result follows  
from
the same commuting square argument that gives Proposition 4.4 and  
Theorem 4.1 of \cite{alr}. Notice that
\cite[Lemma 4.2]{alr} is not needed here because we have factored out  
the action of $\os$.
\end{proof}
\begin{theorem}
  Let $f\mapsto \hat{f}$ denote the Gelfand--Fourier transform of  
$C^*(\kmo)$ onto $C(\invd)$
 obtained from the pairing defined in \proref{dualitypairing}. Denote
 the canonical generators of $C^*(\kmo)$ by $\delta_r$, and let $\{v_a: a\in \ox/\os\}$
be the canonical semigroup of isometries in
 $C(\Omega) \rtimes_\alpha (\ox / \os)$.
Then the maps
\begin{align*}
\teta_r &\mapsto \hat\vartheta_r  := \frac{1}{R(r)}  
\sum_{u\in\os/\stabo{r}}\hat\delta_{ur} & r&\in \kmo\\
\mu_a  & \mapsto v_a &  a&\in \ox/\os,
\end{align*}
determine an isomorphism of
$C^*_r(\pk,\po)$ onto  $C(\Omega) \rtimes_\alpha (\ox / \os)$.
\end{theorem}

\begin{proof}
Part (4) of  \proref{P:Ateta characterization} and the Gelfand--Fourier  
transform give an injective homomorphism of
$A_\theta:=C^*(\{\teta_r: r \in \kmo\})$ into
$C(\invd)$. The range of this homomorphism is (canonically isomorphic to)
$C(\Omega)$, because the conditional expectation  $E_{\osbar}$ is contractive, and thus
transforms the set of generators $\hat\delta_r$, whose linear span is dense in  
$C(\invd)$, into the set of $\osbar$-invariant functions
$\hat\vartheta_r$, whose linear span is dense in~$C(\Omega)$.

Let $\pi_{\teta}$ be the inverse of this isomorphism; the pair
$(\pi_{\teta},\mu)$ is covariant, so the faithfulness criterion
  of \proref{dual+faith}
implies that $\pi_\theta\times \mu$ is an isomorphism.
\end{proof}

\begin{remark}
When $\kk = \Q$,  the unit group is $\os = \{\pm1\}$ and
\cite[Remark 33.b]{bos-con} shows that  the
Hecke algebra of $\po\subset \pk$ is the subalgebra of
the Bost-Connes algebra of fixed points under conjugation by  
$\smatr{0}{-1}$.
  It is also interesting to compare our Hecke C*-algebra to that
  of the almost normal inclusion
  \[
\Gamma_\oo : = \matr{\oo}{1} \quad  \subset \quad   \Gamma_\kk: =  
\matr{\kk}{\kk^*} ,
\]
  considered in \cite{alr}.
Indeed, by \cite[Corollary 2.5]{alr},
$C^*_r(\Gamma_\kk,\Gamma_\oo)$ is canonically isomorphic to
the semigroup crossed product
$C^*(\kmo) \rtimes \ox$, and we can use the embedding of $A_\theta$ in  
$C^*(\kmo) $
to obtain an embedding of crossed products and hence of Hecke algebras.
Since the semigroups in the two crossed products are not the same,
  we need to use the cross section of the semigroup homomorphism $\ox  
\to \ox/\os$
obtained in \lemref{L: section}.
\end{remark}

\begin{proposition}\label{P:relate to ALR}
Let $C^*_r(\Gamma_\kk,\Gamma_\oo)$ be the Hecke C*-algebra of the
almost normal inclusion $\Gamma_\oo : = \smatr{\oo}{1}   \subset
  \smatr{\kk}{\kk^*} =: \Gamma_\kk $,
    denote the generators of $C^*(\kmo)$ by $\delta_r$ and let the  
$v_a$
realize $\ox$ as a semigroup of isometries in  
$C^*_r(\Gamma_\kk,\Gamma_\oo)$.
For each multiplicative cross section $s: \ox/\os \to \ox$ as in  
\lemref{L: section},
the maps
\[
\mu_a \mapsto v_{s(a)} \quad \text{ and } \quad \teta_r \mapsto
\vartheta_r := \frac{1}{|\os/\stabo{r}|} \sum_{u\in\os /\stabo{r}}  
\delta_{ur}
\]
determine an isomorphism of $\heck$ onto a C*-subalgebra $H_s$ of  
$C^*_r(\Gamma_\kk,\Gamma_\oo)$.
The fixed point algebra of $\osbar$ decomposes as
\[
C^*_r(\Gamma_\kk,\Gamma_\oo)^{\osbar} \cong H_s \otimes C^*(\os).
\]
\end{proposition}

\begin{proof}
We have already seen that $\teta_r \mapsto \frac{1}{R(r)} \sum_u
\delta_{ur}$ determines an embedding of $A_\theta$ in $C^*(\kmo)$, and  
it is easy to check that
the elements $ v_{s(a)}$ for $a\in \ox/\os$ satisfy the relations (I)  
and (III),
giving a homomorphism of $\heck$ onto a C*-subalgebra $H_s$ of  
$C^*_r(\Gamma_\kk,\Gamma_\oo)$.
That this is an isomorphism follows from the faithfulness criterion of  
\proref{dual+faith}.

For the last assertion observe that  $\{v_u\colon u\in \os\}$ generates  
a copy of $C^*(\os)$
  that commutes with $H_s$, so it suffices to
prove that  $C^*_r(\Gamma_\kk,\Gamma_\oo)^{\osbar}$ is generated by  
elements of the form
$v_u$ with $u\in \os$ and $v_{s(a)}^* \vartheta_r v_{s(b)}$ with $a,b  
\in \ox/\os$ and $r\in \kmo$. By \lemref{aside} the $\osbar$-average of
$v_a^* \delta_r v_b  $ can be computed using the finite average  
defining $\vartheta_r$,
so $E_{\osbar} (v_a^* \delta_r v_b)  = v_a^* \vartheta_r v_b $.  
Multiplying as needed by
$v_u$ one can replace $a$ and $b$ by $s(a)$ and $s(b)$.
\end{proof}

We point out that the C*-algebra $H_s$ can be thought of as the Hecke  
algebra of a `$+$ construction' associated to the cross section $s$,  
cf. \cite{har-lei,coh}.
The proposition shows that $H_s$ does not depend on $s$ up to  
isomorphism.

\section{A phase transition theorem.}
\label{S:classnumberone}
We consider the action of $\R$ by automorphisms  $\sigma_t$ of
$C(\Omega) \rtimes_\alpha (\ox/\os)$ induced from the dual action $\hat  
\alpha$ via the norm. Specifically, this
action is given by
\begin{gather}
\sigma_t(v_a^* f v_b) = (N_b/N_a)^{it} v_a^* f v_b
\end{gather}
on the monomials spanning the crossed product, where $a,b \in \ox$ and  
$f\in C(\Omega)$.
The existence of such a continuous action $\sigma$ can be verified using
the universal property of the semigroup crossed product or
the presentation of $\heck$, via \proref{dual+faith}.
One simply observes that for each $t\in \R$ the
families $\{\theta_r: r\in \kmo\}$ and $\{N_a^{it} \mu_a: a\in  
\ox/\os\}$
also satisfy the relations and thus induce an automorphism
$\sigma_t$ of $\heck$; the group property is easy to check on  
generators,
and a standard `$\varepsilon/3$-argument' shows that $\sigma$
is continuous.
  It is also possible, and perhaps more interesting,  to implement the  
dynamics spatially
in the Hecke representation on $\ell^2(\po\backslash \pk)$ using
the strongly continuous one-parameter unitary group $(U_t )_{t\in \R}$
on $\ell^2(\po\backslash \pk)$  defined by
  $U_t \xi (\gamma) := (L(\gamma)/R(\gamma))^{-it} \xi(\gamma) $.

  The symmetries are defined by multiplication:
  for $\chi \in W/\osbar$ let  $f_\chi (x) := f(\chi x)$ for $x\in  
\Omega$; then
  \[
  \gamma_\chi(v_a^* f v_b) = v_a^* f_\chi v_b
  \]
  defines a natural strongly continuous action $\gamma$ of $W/\osbar$ on  
$C(\Omega) \rtimes (\ox/\os)$.
Since $\gamma$ clearly commutes with $\sigma$,
we can interpret the automorphisms $\{\gamma_\chi\colon \chi\in  
W/\osbar\}$
as symmetries of the C*-dynamical system $(C(\Omega) \rtimes (\ox/\os),  
\sigma)$.
Through \proref{dual+faith}, the duality pairing allows us to transpose
  the dynamics and the symmetries from
      $C(\Omega) \rtimes (\ox/\os)$ back to the Hecke C*-algebra $\heck$,
where the dynamics is given by
      \[
      \sigma_t(\mu_a^* \teta_r \mu_b) : = (N_b/N_a)^{it} \mu_a^* \teta_r  
\mu_b,
      \]
      and the symmetries are given by
      \[
      \gamma_\chi (\theta_r) = \theta_{\chi r}.
      \]
An explanation is in order to make sure that the product of $\chi  
\osbar \in W/\osbar$
by $(r+\oo)/\os \in (\kmo)^{\os}$ is a well defined class in  
$(\kmo)^{\os}$:
\[
\gamma_\chi (\hat{\teta}_r)(z) = \mathfrac{1}{R(r)} \sum_{u\in  
\os/\stabo{r}} \hat{\delta}_{ur}(\chi z)
= \frac{1}{R(r)}\sum_{u\in\os /\stabo{r}}\chi_{\A (\Q )} (\Tr_{\A(\kk)/\A(\Q)} (ru \chi z )).
\]

We will also need to transpose the canonical dual action $\hat{\alpha}$  
of $\widehat{ \ks/\os}$
on the semigroup crossed product $A_\theta \rtimes (\ox/\os)$ of \proref{actionalpha}
 to the Hecke C*-algebra, via the isomorphism of \thmref{T: universal presentation}.
  We denote this transposed action
also by $\hat\alpha$; it is given by $\hat\alpha_\chi(\mu_a) = \chi(a)  
\mu_a$, and $\hat\alpha_\chi(\theta_r) = \theta_r$, and its fixed point  
algebra is
  $A_\theta :=C^*(\{\theta_r: r\in (\kmo)^{\os}\})$.
We emphasize that the first assertion of the theorem is that the
  equilibrium condition forces full $\hat\alpha$--invariance on the KMS  
states.
  This extra symmetry is a consequence of the stability of equilibrium,  
and
is a stronger property than $\sigma$--invariance alone, since
$\overline{\sigma_\R}$ is strictly smaller than $\hat\alpha_{  
\widehat{\ks/\os}}$
when the norm is not injective on principal integral ideals.

\begin{theorem}\label{T: phasetrans for classone}\label{main}
Suppose $\kk$ is an algebraic number field with class number $h_\kk =  
1$.
Let\/ $(\heck, \sigma)$ be the Hecke C*-dynamical system associated to  
the almost normal inclusion
$\po\subset \pk$, and let \/ $\gamma$ be
the action of\/ $W/\osbar$ as symmetries of this system. Then all  
KMS$_\beta$ states are
  $\hat\alpha$-invariant, and hence determined by their
  values on the generators $\teta_r$ of $A_\theta \subset \heck$.  
Moreover,
\begin{enumerate}
\item[(i)] for each $\beta \in [0,\infty]$ there exists a unique  
$W/\osbar$-invariant
KMS$_\beta$ state $\phi_\beta$ of\/~$\sigma$, given by
\[
\phi_\beta( \teta_{r} ) = N_b^{-\beta} \prod_{p \divides b}
\left(\mathfrac{1 - N_p^{\beta-1}}{1 - N_p\inv}  \right) ,
\]
with $r = a/b$, for $a,b\in\ox$, in reduced form.
\item[(ii)] For $\beta \in [0,1]$ the state $\phi_\beta$ is the unique  
KMS$_\beta$ state for $\sigma$; it
  is a type III hyperfinite factor state when $\beta\neq 0$.

\item[(iii)] For $\beta \in (1,\infty]$, the extreme KMS$_\beta$
states are indexed by $\partial\Omega$,
  the \/$\osbar$-orbits of extreme points in the inverse different:
\begin{itemize}
\item
The extreme ground states are pure and faithful, and are given by
\[
\phi_{\chi,\infty} (\teta_r) = \mathfrac{1}{R(r)} \sum_{u\in  
\os\!/\stabo{r}} \langle r , u \chi \rangle
\]
where $R(r) = |\os/\stabo{r}|$,  the element $\chi \in \partial\invd$  
is a representative
  of an $\osbar$--orbit,  and\/
$\langle\cdot,\cdot \rangle$ indicates the duality pairing of\/ $\kmo$  
and\/ $\invd$
from \proref{dualitypairing};
states corresponding to different $\osbar$--orbits
are mutually inequivalent.

\item For  $\beta \in (1,\infty) $ the extreme KMS$_\beta$ states are  
given by
\[
\phi_{\chi,\beta} (\teta_r) =
\mathfrac{1}{\zeta_\kk(\beta)} \sum_{a\in \ox/\os} N_a^{-\beta}
\big(\mathfrac{1}{R(r)} \sum_{u\in \os/\stabo{r}} \langle ar, u\chi  
\rangle \big),
\]
where $\zeta_\kk$ denotes the Dedekind zeta function of $\kk$.
\end{itemize}
\noindent The state $\phi_{\chi,\beta} $ is quasiequivalent to  
$\phi_{\chi,\infty}$ and the map
$T_\beta: \phi_{\chi,\infty} \mapsto \phi_{\chi,\beta}$
extends to an affine isomorphism of the simplex of KMS$_\infty$ states  
onto the
KMS$_\beta$ states.
The  action of the symmetry group  $W/\osbar$ on the extreme  
KMS$_\beta$ states,
given by $\gamma_w(\phi_{\chi,\beta}) =\phi_{\chi,\beta} \circ \gamma_w  
= \phi_{\chi w, \beta}$,
is free and transitive.

\item[(iv)] The eigenvalue list of the Hamiltonian $H_\chi$ associated  
to $\phi_{\chi,\infty}$ is
$\{\log N_a\colon a\in \ox/\os\}$ for every $\chi$, so the `represented  
partition function'\/
$\operatorname{Tr} e^{\beta H_\chi}$ is independent of\/ $\chi$ and  
equals the Dedekind
zeta function $ \zeta_\kk(\beta)$ of\/ $\kk$.
\end{enumerate}
\end{theorem}

\begin{proof} Since $h_k = 1$,  we have that the principal integral
ideals $\ox/\os$ form a lattice semigroup, and these
act by injective endomorphisms that respect the lattice structure by  
\proref{actionalpha}.
By \cite[Theorem 12]{diri}  KMS$_\beta$ states of $\sigma$ are  
$\hat\alpha$ invariant and
their restrictions to $C(\Omega)$ are characterized by the rescaling  
property
\[
\phi \circ \alpha_a = N_a^{-\beta} \phi.
\]
Using the embedding of $C(\Omega)$ in $C^*(\kmo)$,
and the formula for the conditional expectation
of the action of $W$ on $C^*(\kmo)$ given in \cite[pp.~376]{diri},
one can write the conditional expectation $E_{\gamma}$ explicity:
for $r = a/b$ in reduced form, with $b = \prod_p p^{k_p}$ one has
\[
E_{\gamma} \teta_{a/b}
= \prod_p \left(\frac{N(p)}{N(p) -1} \alpha_{p^{k_p}}(1) -
            \frac{1}{N(p)-1} \alpha_{p^{k_p-1}}(1) \right),
\]
where the product is over a choice
of representatives of the prime ideals in $\oo$, and is independent of  
the choice.
Part (i) follows from this and the rescaling property above, since a  
symmetric KMS$_\beta$ state
must come from a KMS$_\beta$ state on the symmetric system.
The  C*-algebra of the symmetric system, namely,  the range of  
$E_{\gamma}$, is isomorphic to
  $C^*(\ox/\os)$, and its fixed point algebra under $\hat\alpha$ is
  a cartesian product  of sequence spaces (over the nonarchimedean  
valuations $\mko$),
  on which the rescaling condition
is satisfied, for a given  $\beta$, only by the obvious product state.
The KMS$_\beta$ state $\phi$ given in part (i)  is then obtained
by lifting  this product state via $E_{\gamma \times \sigma}$, as in
  \cite[Theorem 34]{diri}.

It is straightforward to characterize the set of $\osbar$--orbits of points in the extreme  
inverse different:
as a subset of the orbit space of the different, it is the  
complement of the (orbits of) nontrivial multiples:
\[
\partial \Omega = \bigcap_{p\in P} (p \invd )^c /\osbar
=  \bigcap_{p\in P} \supp( 1 - \widehat{ \alpha_p(1)}) / \osbar;
\]
so parts (iii)  and (iv) follow  from \cite[Theorem 20 and Corollary  
22]{diri}
and the discussion about the partition function following those results.
We point out that the characterization from \cite{diri} applies to {\em  
ground states},
but since by part (ii) every ground state
is a weak limit of KMS$_\beta$ states as $\beta$ tends to $\infty$,
the distinction of ground states vs.\ KMS$_\infty$ states does not arise.

Part (ii) is proved following the same arguments as in
Neshveyev's ergodicity proof of the rational case \cite{nes}, for which  
we first
need to compute the dilation of the action of $\ox/\os$ on $\Omega$.

Denote by  $\A^0(\kk)^{\osbar}$ the orbit space of the action of  
$\osbar$ on the finite adeles $\A^0(\kk)$. Then
$\A^0(\kk)^{\osbar }$ is a locally compact space containing $\Omega$ as  
a compact subset,
and the action of $\ox/\os$ on $\Omega$ extends to an action of  
$\ks/\os$ on $\A^0(\kk)^{\osbar }$.

\begin{proposition}\label{fullcorner}
The semigroup crossed product $C(\Omega) \rtimes (\ox/\os)$ is  
canonically isomorphic
to the full corner of $C_0(\A^0(\kk)^{\osbar}) \rtimes (\ks/\os)$  
determined by the projection
$1_\Omega \in C_0(\A^0(\kk)^{\osbar }) $.
\end{proposition}
\begin{proof}
By uniqueness of the minimal automorphic extension $\tilde{\alpha}$ of the semigroup  
action $\alpha$ \cite[Theorem 2.1]{minautex},
 it suffices to check that $\{\tilde{\alpha}_q( f)\in  
C_0(\A^0(\kk)^{\osbar }) : q\in \ks/\os, f\in C(\Omega)\}$ is dense in
$C_0(\A^0(\kk)^{\osbar }) $, and this is easy to see from the density  
of $\bigcup_q q\inv \Omega$
in $\A^0(\kk)^{\osbar }$.
\end{proof}

\begin{proposition}
Let  $\mu$ on $\A^0(\kk)^{\osbar }$ be a positive measure, which we  
view, indistinctly,
also as a positive linear functional on $C_0(\A^0(\kk)^{\osbar })$. If  
$\mu$ satisfies
\begin{align}\label{onebeta}
\mu(\Omega) = 1 \qquad \text{ and } \quad q_\star\mu = N_q^\beta \mu,  
\quad \text{where} \quad q_\star\mu(X) = \mu(q\inv X),
\end{align}
then the state $\phi_\mu$ of $\heck$ obtained by restricting $\mu\circ  
E_\gamma$ is a KMS$_\beta$ state.
Conversely, every KMS state arises this way, in fact, the map  
$\mu\mapsto \phi_\mu$ is an affine isomorphism of simplices.
\end{proposition}
\begin{proof}
Direct from the characterization of KMS states given in \cite{diri},  
and the characteristic
property of the minimal dilation/extension; see also \cite{nes}.
\end{proof}

\begin{proposition}
For each $\beta\in (0,1]$ and any measure $\mu$ satisfying  
\eqref{onebeta}, the action
of $\ks/\os$ on $\A^0(\kk)^{\osbar }$ is ergodic.
\end{proposition}
\begin{proof}
The proof parallels that of \cite[Proposition]{nes},
using \proref{fullcorner} in place of \cite[Corollary 3.2]{minautex}.
One has to implement the following changes in Neshveyev's proof:
\begin{eqnarray*}
\mathcal P &\longrightarrow& \mko  \\
\A^0 &\longrightarrow& \A^0(\kk)^{\osbar }\\
R:= \prod \Z_p &\longrightarrow& \Omega:=  (\prod \invdv) /\osbar\\
\text{(origins of orbits) }R^* &\longrightarrow&  \partial \Omega \\
\text{(symmetries) }W: = \prod \Z_p^* &\longrightarrow&  W/\osbar =  
(\prod \oo_v^*) /\osbar.
\end{eqnarray*}
Notice that the last two items in the list coincide for $\kk = \Q$
(on the left column), but  
have
different versions for a general number field (on the right column). The crucial ingredients  
are the convergence of the Dirichlet $L$--series $L(\omega,\beta)$, and  
the divergence of the Dedekind zeta function as $\beta \to 1^+$,   
\cite[Theorem 5, Ch. VIII, p. 161]{lan}.
\end{proof}

Since for $\beta \in (0,1]$ the measure corresponding to the
symmetric KMS$_\beta$ state $\phi_\beta$ is ergodic,
and KMS$_\beta$ states form a convex set, the
uniqueness  assertion  in part (i) follows.
This finishes the proof of \thmref{main}.
\end{proof}

\begin{remark}
With respect to part (iv) of \thmref{main}, we point out that by  
\cite[Theorem 1]{ole-ped},
the dynamics $ \sigma$ is not approximately inner, so there
is no intrinsic Hamiltonian for this system. But the Liouville operators
associated to the extreme KMS$_\infty$ states 
(i.e.\ the `Hamiltonian operators' in the associated GNS representations)
have all the same spectrum,
so it makes sense to talk about the spectrum of the Hamiltonian and
to associate a partition function to the system; see
\cite{win-ten} for a discussion of the spectra of Liouville operators.
\end{remark}

\section{Class field theory considerations}
\label{S:Classfields}
Let $H$ be the Hilbert class field of $\kk$,
and let $H^+$ be the maximal abelian extension of $\kk$ in which no  
finite prime is ramified.
For example, for $\Q$, we have $H=H^+=\Q$,
but $\kk=\Q(\sqrt{3})$ has $H=\kk$ and $H^+=\kk(i)$.

The symmetry group $W/\osbar$ has a Galois group interpretation as  
$\mathcal G(\kk^{ab},H^+)$,
the Galois group of the maximal abelian extension of $\kk$ over $H^+$.
We embed $W$ into $\mathbb A_\kk^*$ by  
$(a_v)_{v\in\mko}\longmapsto(1,\dots,1,a_v)_v$,
where we put $1$ in each archimedean component, cf. \cite[Proposition  
8.5]{har-lei}.

\begin{proposition}\label{exactsequence}
Let\/ $\sigma_1,\dots,\sigma_r$ be the embeddings of\/ $\kk$ in $\R$,
and let\/ $\sgn(\os)$ be the subgroup of\/ $\left<-1\right>^r$  
generated by the elements
$(\sgn(\sigma_iu))_{i=1,\dots,r}$ (\/$u\in\os$).
Then the sequence
\[
1\longrightarrow W/\osbar\longrightarrow\mathbb A_\kk^*/  
\overline{\kk^* \kk^1_\infty}
\longrightarrow\Cl(\kk)\times\left<-1\right>^r/\sgn(\os)\longrightarrow1
\]
is exact. Via the Artin map, we obtain the canonical isomorphism
\begin{gather}
W/\osbar \quad \cong \quad  \mathcal G(\kk^{ab}:H^+).
\end{gather}
\end{proposition}

\begin{proof}
The exactness on the left follows from the proof of \lemref{aside}.
The group of finite ideles modulo $W$ is isomorphic to the group of  
ideals,
and modulo~$\kk^*$ we may change an idele to one of a fixed set of  
representatives of $\Cl(\kk)$.
Then we can still change the infinite part of the idele by an element  
of $\os\kk^1_\infty$,
which means that we can change the signs at the real places by an  
element of $\sgn(\os)$.

The isomorphism $\mathbb A_\kk^*/ \overline{\kk^* \kk^1_\infty} \cong  
\mathcal G(\kk^{ab}: \kk)$
is from \cite[p.~272]{wei74}.
Since the group on the right is a finite factor group that contains no  
nontrivial unit at a finite
  place,
it corresponds to a finite extension of $\kk$ in which no finite place  
is ramified.
Since~$\os$ acts trivially,
all other places are allowed to ramify.
Hence this group is the Galois group of $H^+$ over $\kk$,
and, via the Artin map, $W/\osbar$ is isomorphic to the Galois group of  
$\kk^{ab}$ over $H^+$.
\end{proof}

\begin{corollary}\label{C: pureimag-classone}
If $\kk$ is a number field of class number one with no real embeddings,
then $W/\osbar$ is isomorphic to $G(\kk^{ab}: \kk)$.
Via this isomorphism, 
the extreme KMS$_\beta$ states of the system $(\heck,\sigma)$ for $\beta > 1$
are indexed by the complex embeddings of the maximal abelian extension  
$\kk^{ab}$ of\/ $\kk$.
\end{corollary}
\begin{proof}
  $H=\kk$ because $h_\kk=1$ and
$H^+ = H$ because $\kk$ has no real embeddings, so the isomorphism
of $W/\osbar$  to $G(\kk^{ab}: \kk)$ follows at once from \proref{exactsequence}. 
By \thmref{main} part (iii), the action of $W/\osbar$ is free and transitive on 
extreme KMS$_\beta$ states, so the second assertion follows.
\end{proof}

\begin{corollary}
If the class number of $\kk$ is still one,
but there are nontrivial unramified extensions (at the finite places),
then extreme KMS$_\beta$ states  for $\beta > 1$
are indexed by the complex embeddings of the maximal abelian extension,
modulo the complex conjugations over each real embedding of $\kk$.
\end{corollary}

\begin{proof}
The result follows from the fact that  
$\Gal(\kk^{ab}:H^+) \cong \Gal(\kk^{ab}: \kk) / (\{\pm 1\} ^{r_1} )$,
where $r_1$ denotes the number of real embeddings of $\kk$.
\end{proof}

One can use the indexing of \corref{C: pureimag-classone}
to obtain an action of $G(\kk^{ab}: \kk)$
on extreme KMS$_\infty$ states, essentially coming from the 
`geometric' symmetries of the dynamical system. 

It is also possible to define
another action of $G(\kk^{ab}: \kk)$
on extreme KMS$_\infty$ states, coming from `arithmetic' symmetries.
First define the {\em arithmetic Hecke algebra} $\kk (\pk,\po)$ to be the  
algebra over $\kk$
generated by the $\theta_r$ and the $\mu_a$. Evaluation of an extreme
  KMS$_\infty$ state on a $\theta_r$
corresponds to evaluation of $\hat\theta_r$ on the
$W/\osbar$--orbit of a point in the  extreme inverse different. Since  
$\hat\theta_r$
is a $\Q$--linear combination of characters, it follows that
  the image of $\kk (\pk,\po)$ under
  an extreme KMS$_\infty$ state is
  contained in the real subfield of the maximal {\em cyclotomic}  
extension of $\kk$. Thus, there is an action of $G(\kk^{ab}: \kk)$
on the values of extreme  KMS$_\infty$ states.
Since in general this is a proper subfield of $\kk^{ab}$,
 the arithmetic Hecke algebra cannot support 
 a canonical action in which $G(\kk^{ab}: \kk)$ 
 acts freely and transitively by composition with extreme  
KMS$_\infty$ states.
  In fact, one can compute both actions and see how they differ.

\begin{theorem}\label{badnews}
 The geometric action of $W/\osbar$ on extreme KMS$_\infty$ states 
from \corref {C: pureimag-classone} coincides with the Galois  
action of $W/\osbar$ on their values, $\phi_{\chi,\infty} (\teta_r)$, as given by  
class field theory, if and only if $\kk =\Q$.
\end{theorem}

\begin{proof}
Let $\phi_{\chi,\infty}$ be an extreme KMS$_\infty$ state, and let $j$  
be an  idele.
The action of~$j$, viewed as a symmetry of $\Omega = \invd /\osbar$,
on $\phi_{\chi,\infty}$ is given by
\[
\phi_{\chi,\infty} (\teta_r) \mapsto  \phi_{j \chi,\infty} (\teta_r) =
\frac{1}{R(r)}\sum_{u\in\os /\stabo{r}}\chi_{\A (\Q )} (\Tr_{\A(\kk)/\A(\Q)} (j ru\chi  
)).
\]

Next we compute the action of the idele $j$,
viewed now as a Galois element acting on the values of the extreme KMS  
states
via the Artin map as in \proref{exactsequence}.
We note first that the complex number $\phi_{\chi,\infty} (\teta_r)$ is  
a linear combination
with integer coefficients of character values of the group  
$\kmo$.
Thus this combination of roots of unity is in the maximal cyclotomic  
extension of $\Q$.
By class field theory, the Galois action of an idele $j\in W$, when  
restricted to $\Q^{cycl}$,
coincides with the action of $N(j)$,
where $N(j)$ is the norm of $j$ to the rational ideles
(see~\cite[Corollary~1, p.~246]{wei74}).
Thus the Galois action of $j$ on values of KMS$_\infty$ states is given  
by
\[
\phi_{\chi,\infty} (\teta_r)
\mapsto\frac{1}{R(r)}\sum_{u\in\os /\stabo{r}}\chi_{\A (\Q )} (\Tr_{\A(\kk)/\A(\Q)}
(N(j) ru\chi)),
\]
where we have used the $\A (\Q )$-linearity of the trace to replace
$N(j)\Tr_{\A(\kk)/\A(\Q)} (ru\chi )$ by $\Tr_{\A(\kk)/\A(\Q)} (N(j)ru\chi )$.

If we now take $j$ to be a rational idele, then $N(j)=j^d$, where  
$d=\deg (\kk /\Q )$,
and it follows that the two actions are different unless $d = 1$,
 i.e.\ unless $\kk =\Q$.
\end{proof}

In order to relate our results to those of \cite{har-lei}, we
recall that the Hecke algebras of \cite{har-lei} depend in general on  
two choices:
first there is a localization process that substitutes the ring of  
integers by a principal ring,
and then there is a multiplicative cross section from the principal  
ideals into the ring.
Thus the resulting  almost normal subgroup is not canonical,  and
one should not expect a general direct relation with the Hecke algebras
considered here and in \cite{alr}, which are canonical.

However, if we assume $h_\kk = 1$ and fix a subsemigroup $S$ of $\ox$  
containing
  a representative for each ideal of $\oo$ (see \lemref{L: section}),
then the Hecke C*-algebra $C^*(\Gamma_S,\Gamma_\oo)$ of the almost  
normal inclusion
\[
\Gamma_\oo := \matr{\oo}{1} \subset \matr{\kk}{SS\inv} = : \Gamma_S
\]
considered in \cite{har-lei} has an action of $W = \prod_{v\in \mko}  
\ov^*$ by
\cite[Section 4]{har-lei}, and one obtains the following
generalization of part~3 of \cite[Proposition 8.5]{har-lei}.
\begin{proposition}
Suppose $h_\kk =1$, and let $S$ be the range of
  a multiplicative cross section for $\ox \to \ox/\os$. Then
\[
C^*(\Gamma_S,\Gamma_\oo)^{\osbar} \cong \heck
\]
\end{proposition}
\begin{proof}
Embed $C^*(\Gamma_S,\Gamma_\oo)$ in $C^*(\Gamma_\kk,\Gamma_\oo)$ as  
indicated in
\cite[Remark 5.4]{alr} . This shows that $C^*(\Gamma_S,\Gamma_\oo)$
does not depend on the specific cross section $S$ up to isomorphism.
The result now follows from \proref{P:relate to ALR}.
\end{proof}

\begin{remark}\label{lastremark}
The particular case of
\corref{C: pureimag-classone} corresponding to the nine
quadratic imaginary fields of class number one,
 i.e.\ $\kk = \Q[\sqrt{-d}]$
 for $d = 1,\  2, \  3  ,$ $7,$ $11,$ $19,$ $43,$  
$67,$ $163$,
is already implicit in \cite{har-lei}.
Indeed, for such fields $\os$ is a finite group, and
the almost normality of the inclusion $\po \subset \pk$ is
straightforward, and does not require the general considerations leading
to \lemref{L:almostnormal}. As indicated by Harari and Leichtnam,
the argument of \cite[Remark 33.b)]{bos-con} can be used to show that  
$\heck$ is the
fixed point algebra of their Hecke C*-algebra $C^*(\Gamma_S,\Gamma_\oo)$
under the action of $\os$, see also \cite[Example 2.9]{LaLa}.
   That the extreme KMS$_\beta$ states of $\heck$ for $\beta > 1$
are indexed by  $\mathcal G(\kk^{ab}: \kk)$ then
follows from \cite[Theorem~0.2 and Proposition 8.5]{har-lei}.
This observation and the final remarks \cite[pp.~241--242]{har-lei}
concerning the Artin map, were a strong motivation for the present work,
by reinforcing our belief that the almost normal inclusion of full  
``$ax+b$'' groups considered
here would lead to a C*-dynamical system with group
of symmetries isomorphic to the right Galois groups, and that
an appropriate compactification of $\os$, like the one given in  
\lemref{aside},
was the key to understand the corresponding Hecke C*-algebra.
\end{remark}

 \end{document}